\documentclass[12pt,a4paper]{article}
\usepackage[utf8]{inputenc}
\usepackage{pgf,tikz}
\usepackage{mathrsfs}
\usetikzlibrary{arrows}
\usepackage{yfonts}

\usepackage[cyr]{aeguill}
\usepackage[francais]{babel}

\DeclareMathSymbol{A}{\mathalpha}{operators}{`A}
\DeclareMathSymbol{B}{\mathalpha}{operators}{`B}
\DeclareMathSymbol{C}{\mathalpha}{operators}{`C}
\DeclareMathSymbol{D}{\mathalpha}{operators}{`D}
\DeclareMathSymbol{E}{\mathalpha}{operators}{`E}
\DeclareMathSymbol{F}{\mathalpha}{operators}{`F}
\DeclareMathSymbol{G}{\mathalpha}{operators}{`G}
\DeclareMathSymbol{H}{\mathalpha}{operators}{`H}
\DeclareMathSymbol{I}{\mathalpha}{operators}{`I}
\DeclareMathSymbol{J}{\mathalpha}{operators}{`J}
\DeclareMathSymbol{K}{\mathalpha}{operators}{`K}
\DeclareMathSymbol{L}{\mathalpha}{operators}{`L}
\DeclareMathSymbol{M}{\mathalpha}{operators}{`M}
\DeclareMathSymbol{N}{\mathalpha}{operators}{`N}
\DeclareMathSymbol{O}{\mathalpha}{operators}{`O}
\DeclareMathSymbol{P}{\mathalpha}{operators}{`P}
\DeclareMathSymbol{Q}{\mathalpha}{operators}{`Q}
\DeclareMathSymbol{R}{\mathalpha}{operators}{`R}
\DeclareMathSymbol{S}{\mathalpha}{operators}{`S}
\DeclareMathSymbol{T}{\mathalpha}{operators}{`T}
\DeclareMathSymbol{U}{\mathalpha}{operators}{`U}
\DeclareMathSymbol{V}{\mathalpha}{operators}{`V}
\DeclareMathSymbol{W}{\mathalpha}{operators}{`W}
\DeclareMathSymbol{X}{\mathalpha}{operators}{`X}
\DeclareMathSymbol{Y}{\mathalpha}{operators}{`Y}
\DeclareMathSymbol{Z}{\mathalpha}{operators}{`Z}
\DeclareMathSymbol{e}{\mathalpha}{operators}{`e}

\newcommand{\mathbb}{\mathbf}

\newcommand{\KK}{\mathbb{K}}
\newcommand{\LL}{\mathbb{L}}

\newcommand{\RR}{\mathbb{R}}

\gdef\thechapter{\@Roman\c@chapter}%

\newtheorem{theoreme}{\noindent {\textbf{Th\'eor\`eme}}}[section]
\newtheorem{proposition}{\noindent \textbf{Proposition}}[section]

\newtheorem{definition}{\noindent \textbf{D\'efinition}}[section]

\newenvironment{demonstration}{\begin{trivlist}\item[]{\textbf{D\'emonstration.}}~---}%
{\nolinebreak  \end{trivlist}}

\newcounter{exos}

\newcommand{\cC}{\mathcal{C}}
\newcommand{\cD}{\mathcal{D}}
\newcommand{\voir}{{\emph{voir}}}
\newcommand{\vs}{\vspace*{3mm}}

\title{La notion d'involution dans le {\it Brouillon Project} de Girard Desargues}
\author{
par Marie Anglade et Jean-Yves Briend}
\date{\today}
\begin{document}
\maketitle

\begin{abstract} Nous tentons dans cet article de proposer une thèse cohérente concernant la formation de la notion d'involution dans le {\it Brouillon Project} de Desargues. Pour cela, nous donnons une analyse détaillée des dix premières pages dudit {\it Brouillon,} comprenant les développements de cas particuliers qui aident à comprendre l'intention de Desargues. Nous mettons cette analyse en regard de la lecture qu'en fait Jean de Beaugrand et que l'on trouve dans les {\it Advis Charitables.} 
\end{abstract}

\vspace*{1cm}

\begin{center}
{\Large\bf Introduction}
\end{center}

L'{\oe}uvre mathématique de Girard Desargues a été l'objet de nombreux travaux. Nous pouvons citer la première édition de ses {\oe}uvres donnée par Noël-Germinal Poudra en 1864, celle de René Taton en 1951 ou encore leur édition et traduction en anglais par Judith Field et Jeremy Gray en 1987 ({\it voir} respectivement \cite{poudra1,poudra2}, \cite{taton} et \cite{fieldgray}). Les études du contenu mathématique du {\it Brouillon Project} se concentrent le plus souvent sur la partie du texte portant sur les coniques. Nous pouvons citer par exemple l'article de Jan Hogendijk \cite{hogendijk} qui analyse les apports de Desargues au regard de l'héritage d'Apollonius. Il est également étudié dans ses rapports à la théorie et à l'histoire de la perspective, entre autres par Kirsti Andersen, James Elkins ou Jeanne Peiffer ({\it voir} respectivement \cite{andersen}, \cite{elkins} et \cite{peiffer}). La notion d'involution a beaucoup moins attiré l'attention des chercheurs et à notre connaissance l'article de Frédérique Lenger \cite{lenger}, paru en 1950, est le seul à porter spécifiquement sur ce sujet. Il faut cependant citer Max Zacharias dans \cite{zacharias}, Jean-Pierre le Goff qui traite assez longuement de la notion dans \cite{legoff} et y propose une thèse sur l'heuristique ayant mené Desargues à l'involution, ainsi que Laura Catastini et Franco Ghione dans \cite{catastinighione}, qui en présentent un aperçu assez large au public mathématique italien. Nous voudrions dans ce texte combler ce qui nous semble être une lacune pour la compréhension du rôle de Desargues dans l'histoire de la géométrie en analysant le plus complètement possible la partie du {\it Brouillon Project} qui a trait à l'involution.

\begin{figure}[!h]
\centering
\definecolor{uuuuuu}{rgb}{0.26666666666666666,0.26666666666666666,0.26666666666666666}
\definecolor{xdxdff}{rgb}{0.49019607843137253,0.49019607843137253,1.}
\definecolor{qqqqff}{rgb}{0.,0.,1.}
\begin{tikzpicture}[line cap=round,line join=round,>=triangle 45,x=0.8583690987124469cm,y=0.8583690987124469cm]
\clip(-1.08,1.78) rectangle (17.56,12.38);
\draw [rotate around={31.458346614308038:(9.78,7.45)}] (9.78,7.45) ellipse (4.57642608655354cm and 3.231845753094741cm);
\draw [domain=-1.08:17.56] plot(\x,{(--49.1748405018268--5.4790271290221995*\x)/13.875145648975273});
\draw [domain=-1.08:17.56] plot(\x,{(--27.662904474155628-0.6391262300794445*\x)/7.034094168417636});
\draw [domain=-1.08:17.56] plot(\x,{(--17.12761631830195--5.697129477714695*\x)/5.6179585234387845});
\draw [domain=-1.08:17.56] plot(\x,{(--54.20001804359493-6.336255707794139*\x)/1.4161356449788514});
\begin{scriptsize}
\draw[color=black] (5.36,7.6) node {$\cC$};
\draw [fill=qqqqff] (0.8,3.86) circle (2.5pt);
\draw[color=qqqqff] (0.78,4.28) node {$c$};
\draw [fill=xdxdff] (14.675145648975274,9.339027129022199) circle (2.5pt);
\draw[color=xdxdff] (14.86,9.8) node {$h$};
\draw[color=black] (-0.58,3.58) node {$\cD$};
\draw [fill=uuuuuu] (4.902306542353444,5.479921650248029) circle (1.5pt);
\draw[color=uuuuuu] (5.04,5.86) node {$b$};
\draw [fill=uuuuuu] (6.417958523438784,9.557129477714694) circle (1.5pt);
\draw[color=uuuuuu] (6.08,9.76) node {$p$};
\draw [fill=uuuuuu] (7.834094168417636,3.2208737699205554) circle (1.5pt);
\draw[color=uuuuuu] (8.06,3.54) node {$q$};
\draw [fill=uuuuuu] (7.1323879453122006,6.360537740043187) circle (1.5pt);
\draw[color=uuuuuu] (7.28,6.8) node {$g$};
\end{scriptsize}
\end{tikzpicture}
\caption{Construction d'une division harmonique à partir d'une conique.}\label{Polaire}
\end{figure}
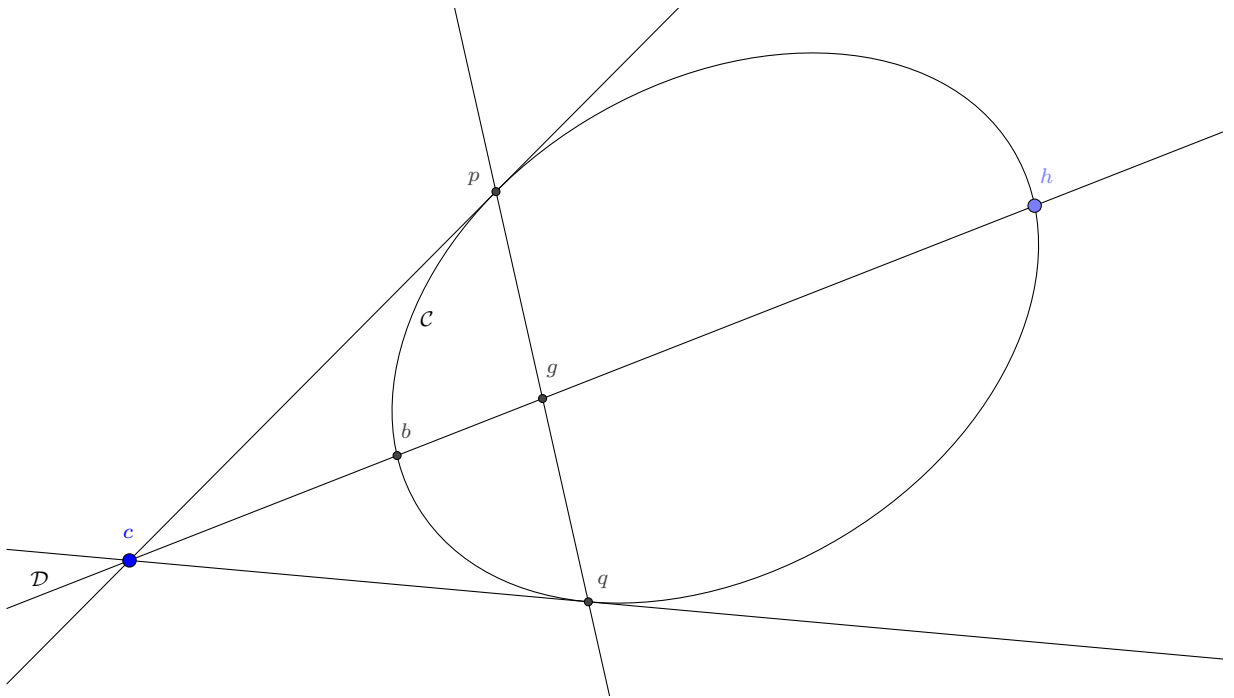

Soient $\cC$ une conique et $\cD$ une droite coupant la conique aux deux points $b$ et $h$. Prenons un point $c$ sur $\cD$ et traçons depuis $c$ les deux tangentes à la conique $\cC$, qui sont tangentes à celle-ci en les points $p$ et $q$ respectivement. La droite $pq$ coupe la droite $\cD$ en un point $g$ ({\it voir} la figure \ref{Polaire}). La théorie apollonienne des coniques\footnote{{\it Voir} par exemple la proposition 37 du livre III des {\it Coniques} d'Apollonius, p. 351 des commentaires de Roshdi Rashed dans \cite{apollonius-rashed}} permet de conclure que la droite $\cD$ est harmoniquement coupée par les quatre points $b,h,c,g$ ou encore que la double paire $(b,h;c,g)$ en est une \emph{division harmonique.} Étant donnée la conique $\cC$ et la droite $\cD$ qui coupe $\cC$ en $b,h$, on dispose, en faisant varier le point $c$ sur $\cD$, d'une foule de divisions harmoniques $b,h,c,g$ dont $b,h$ est l'un des couples. 

On peut se poser alors deux questions. La première est de savoir ce qui lie entre elles ces différentes divisions harmoniques. La seconde est de savoir si, étant donnée une telle collection de divisions harmoniques sur une droite $\cD$, on peut reconstruire une conique les produisant selon la méthode décrite ci-dessus. La réponse à ces deux questions est à la base du renouveau de la théorie des coniques au début du dix-neuvième siècle entre les mains d'Augustin Fresnel, Jean-Victor Poncelet ou Jakob Steiner. 

La première question pourrait se reformuler comme suit. Étant donnée une division harmonique $b,h,c,g$, prenons deux autres points $c',g'$ de sorte que $b,h,c',g'$ soit encore une division harmonique. Peut-on dire qu'alors $c,g,c',g'$ est aussi une division harmonique? La réponse est non, et cela tient au fait qu'il existe, à partir de la division $b,h,c,g$, deux manières de construire un autre couple $c',g'$ en division harmonique avec $b,h$. La relation de division harmonique n'est donc pas {\it agrégative,} au sens où l'est, par exemple, celle d'alignement. Précisons ce que nous entendons par là~:~étant donnés trois points alignés $b,h,c$, si $g$ est un quatrième point aligné avec deux de ces trois points, il l'est avec deux \textit{quelconques} parmi ces trois. Pour  mettre sous un même toit toutes les divisions harmoniques construites sur une droite comme nous l'envisageons ci-dessus, il faut une notion qui soit plus forte que celle de la simple division harmonique, en levant l'ambiguïté mentionnée plus haut concernant la construction d'un tiers couple $c',g'$ tel que la droite soit harmoniquement coupée par $b,h,c',g'$. La notion arguésienne d'{\it involution} s'avère remplir exactement ce contrat. Le point clé est le suivant~:~si l'on considère la notion de division harmonique comme une relation entre deux {\it couples} de points, elle est insuffisante car non agrégative, comme l'est d'ailleurs la relation d'alignement pour deux points. En passant à la notion d'involution de {\it trois} couples de points, nous retrouvons une notion agrégative ayant la même puissance générative que celle d'alignement de trois points. Elle permettra, plus d'un siècle après sa mise au point par Desargues, de répondre efficacement à la seconde question énoncée ci-dessus, menant à la théorie de la dualité, des homographies et des paramétrages rationnels des coniques ({\voir} le livre \cite{chasles} de Michel Chasles), même si la réponse se trouve déjà en germe dans le {\it Brouillon Project,} comme en attestent par exemple les lignes 24 et 25 de la page 29 du \textit{Brouillon.}

Nous allons dans cet article montrer que c'est sans doute par une démarche analogue à celle que nous venons de décrire que Desargues a mis au jour la notion d'involution. Pour étayer notre propos, nous nous baserons sur une analyse détaillée des dix premières pages du {\it Brouillon Project} de 1639, et plus particulièrement sur le long développement qu'il y fait des cas particuliers de cette notion où apparaissent ce qu'il appelle des {\it n{\oe}uds moyens.} Nous montrerons en outre, au travers de l'exemple de Jean de Beaugrand, comment les lecteurs contemporains, hormis Blaise Pascal, sont passés à côté de la puissance de cette notion, en restant trop fidèles aux conceptions classiques d'Apollonius et de Pappus.


\begin{center}
{\Large\bf Partie I

\vs
Les dix premières pages du {\it Brouillon--Project}}
\end{center}

Nous tentons dans cette partie de présenter  le contenu des dix premières pages du {\it Brouillon project d'une atteinte aux evenemens des rencontres du Cone avec un Plan,} de Girard Desargues, en un langage compréhensible par le lecteur contemporain. Nous nous sommes basés sur l'unique exemplaire connu, numérisé par la Bibliothèque nationale de France et conservé au département Réserve des livres rares sous la référence RESM-V-276. Le texte de Desargues est tout à la fois relativement facile à lire, si l'on ne s'attache pas aux détails (on «~voit bien~» en quelque sorte où il veut en venir) et plutôt obscur lorsque l'on cherche à en extraire précisément la démarche de l'auteur. Comme son titre l'indique, c'est un brouillon et Desargues répète en plusieurs occurences qu'il mériterait d'être repris, simplifié et amélioré de bien des manières. Cet aspect inachevé lui est reproché notamment par Jean de Beaugrand dans la lettre ouvrant les {\it Advis Charitables} que nous étudions dans la deuxième partie de cet article. 

Nous allons ici essayer de démêler l'écheveau et pour cela nous emploierons sans vergogne les notations algébriques, que Desargues, malgré l'exhortation de Descartes\footnote{{\it cf.} la lettre de Descartes à Desargues du 19 juin 1639, dans \cite{taton}, p. 185.}, s'est refusé à employer. Pour chaque notion, énoncé ou preuve dont nous traiterons, nous essaierons autant que possible de renvoyer à l'endroit de l'original qui y correspond.

Cette partie est organisée comme suit. Nous commençons par rappeler le vocabulaire introduit par Desargues concernant les droites concourantes et parallèles, où il unifie  ces deux notions par l'introduction de points à l'infini. Nous entrons alors dans le vif du sujet, à savoir l'involution, en séparant différents aspects de cette notion~:~combinatoire, métrique ou affine, puis projective. Cette présentation est bien évidemment très différente de celle que fait Desargues, mais c'est à ce prix que l'on peut comprendre sa démarche, face à un original parfois très confus, voire déroutant. Nous analysons à la suite les notions de n{\oe}uds moyens simples et doubles, qui sont au c{\oe}ur de la machinerie arguésienne. Nous finissons ce chapitre par une présentation faite entièrement en un langage mathématique moderne, celui de la géométrie projective, ce qui sera l'occasion de montrer que certaines des difficultés sur lesquelles Desargues bute sont liées au fait qu'il n'a pas à sa disposition les outils conceptuels permettant de distinguer, suivant les cas qu'il considère, à quel type de structure il fait appel (euclidienne, affine ou projective). Nous renvoyons le lecteur à l'appendice pour le détail des outils euclidiens employés par Desargues dans son traitement de l'involution.

Avant de nous plonger dans la lecture détaillée du {\it Brouillon Project} et sans définir les termes qui seront étudiés en détail plus loin, résumons brièvement le cheminement de Desargues dans ces dix pages. Il commence par unifier les notions de droites concourantes et parallèles en introduisant ce qu'il appelle une \emph{ordonnance} de droites, dont le {\it but} ou point de concours est soit à distance finie soit à distance infinie. Il introduit ensuite sa terminologie arboricole, peut-être comme un écho du vignoble de Château  Grillet, propriété de la famille Desargues~:~brins, rameaux, tronc etc. Il définit la notion d'{\it arbre,} comme une configuration particulière de trois {\it couples} de points et d'une {\it souche} alignés sur une droite qu'il nomme {\it tronc}, définie par des égalités de rectangles (ou d'aires de rectangles, dirions-nous aujourd'hui).  Il montre\footnote{Mais ne démontre pas, {\it voir} plus loin.} alors, par des manipulations classiques de la géométrie euclidienne sur les proportions, que cette notion est équivalente à celle d'{\it involution\footnote{Notons que ce terme, comme les précédents, est employé en botanique, désignant l'état de feuilles s'enroulant sur elles mêmes.},} où ne subsistent que les trois couples de points et d'où la souche a disparu. Desargues entre alors dans plusieurs pages touffues (les pages 5 à 10 du {\it Brouillon Project}) où il examine certains cas particuliers dégénérés de la configuration involutive~:~ceux où certains points se confondent, donnant ce qu'il appelle des {\it n{\oe}uds moyens.} Deux cas se présentent à lui, selon que ce sont deux points d'un même couple qui se confondent (les n{\oe}uds moyens {\it doubles}) ou que ce sont deux couples de points distincts qui se confondent (les n{\oe}uds moyens {\it simples}). C'est ici qu'il découvre une situation de division harmonique qu'il peut interpréter de {\it deux manières différentes} comme une involution. Notons au passage qu'il clarifie complètement la division harmonique\footnote{Terme que Desargues n'emploie pas dans le \textit{Brouillon.}} comme étant une {\it relation} entre deux paires de points, source qui pourrait avoir inspiré Leibniz pour sa géométrie relationnelle dans la {\it Characteristica Geometrica} ({\it voir} l'article \cite{debuiche} de Valérie Debuiche). Il conclut ce long développement au milieu de la page 10 par une phrase cruciale ({\it cf. infra}), montrant que la notion d'involution est bien {\it agrégative} au sens où nous l'avons définie dans l'introduction.

\section{Droites et plans, ordonnances et buts}
Dans les deux premières pages du {\it Brouillon,} Desargues donne quelques précisions sur ce qu'il entend par droite, plan et {\it ordonnance} de ces objets. Le cadre implicite dans lequel il se place est celui de la géométrie euclidienne, dans l'espace ou dans un plan. Une ligne droite est entendue comme étant étendue infiniment de part et d'autre\footnote{p.1, l.11}. De même un plan est entendu comme étant infiniment étendu de toutes parts\footnote{p.1, l.30}. La notion importante introduite ici est celle d'ordonnance~:~
\begin{definition} Des droites sont dites être d'une \emph{même ordonnance} si elles sont ou bien concourantes ou bien parallèles\footnote{p.1, l.15}. De même des plans sont dits être d'une même ordonnance s'ils sont concourants ou parallèles\footnote{p.1, l.33}.
\end{definition}
Dit autrement, soit $(D_i)_{i\in I}$ une famille de droites. Les droites de la famille sont de même ordonnance s'il existe un point commun à toutes ces droites ou si toutes les droites de la famille sont parallèles entre elles. On peut, comme Desargues\footnote{p.1, l.19}, dire que la famille de ces droites \emph{forme une ordonnance.} L'un des principes guidant Desargues semble être d'adopter autant qu'il est possible un point de vue unificateur et cela commence ici avec la notion de {\it but} d'une ordonnance~:~
\begin{definition} Le \emph{but} d'une famille de droites d'une même ordonnance est
\begin{itemize}
\item leur point de concourance si ces droites sont concourantes;
\item le point à l'infini en chacune d'elles si elles sont parallèles.
\end{itemize}
Le \emph{but} d'une famille de plans d'une même ordonnance est
\begin{itemize}
\item leur droite de concourance si ces plans sont concourants;
\item la droite à l'infini en chacun d'eux s'ils sont parallèles.
\end{itemize}
\end{definition}
Les termes de Desargues ne sont pas tout à fait les mêmes. Tout d'abord, il utilise le mot \emph{endroit} pour définir la notion de but~:~le but d'une ordonnance de droite, c'est~«~l'endroit auquel on \emph{conçoit}  que les droites de l'ordonnance tendent toutes\footnote{p.1, l. 19}~»~et ce, indépendamment de la nature (concourance ou parallélisme) de l'ordonnance.  Il précise que dans le cas d'une ordonnance de droites parallèles, cet endroit est à distance infinie, d'une part et d'autre, en chacune d'elles\footnote{p.1, l.22} ce que l'on peut interpréter de la manière suivante~:~chaque droite possède un «~point à l'infini~» et ce point est l'endroit vers lequel on tend, que l'on s'éloigne d'un côté ou de l'autre dans la droite en question. Suivant son humeur unificatrice, il énonce que deux droites d'un même plan sont toujours d'une même ordonnance, ce qui en termes modernes n'est autre que le fait que deux droites d'un même plan projectif sont toujours concourantes\footnote{p.1, l.28}. Il est à noter que le choix du mot {\it ordonnance} pourrait avoir été induit par le terme d'{\it ordonnées} utilisé dans la théorie apollonienne des coniques~:~des ordonnées forment une ordonnance dont le but, chez Apollonius\footnote{\textit{voir} les définitions du livre I des \textit{Coniques} d'Apollonius, par exemple dans \cite{apollonius-I}, p. 254.}, est toujours à l'infini. Il reprendra d'ailleurs le terme d'ordonnées dans son développement sur les coniques, p. 15 du \textit{Brouillon}. C'est un premier indice du fait que Desargues n'a pas, contrairement à ce qu'ont pu écrire ses détracteurs, choisi ses appellations nouvelles avec légèreté.

Concernant les plans, il n'utilise pas la terminologie de «~droite à l'infini~» pour parler du but d'une ordonnance de plans parallèles, ne s'étendant pas vraiment ici sur la nature de ce but. Il précise cependant les choses un peu plus loin\footnote{p.1 l. 47 à p.2, l.15}~:~fixons un point et faisons tourner une droite passant par ce point sur toute sa longueur. Un point toujours situé à même distance du point fixe va tracer un cercle (la ligne «~courbée en pleine rondeur~»). Dit autrement, des droites d'une même ordonnance dont le but est à distance finie déterminent des cercles. Prenons maintenant une ordonnance de droites parallèles entre elles. On peut alors concevoir de faire se mouvoir une des droites de l'ordonnance sur toute sa longueur, en laissant le point à l'infini fixe. Cette droite mouvante va alors parcourir une famille de droites parallèles et chaque point de cette droite va tracer à son tour une droite perpendiculaire à toutes les droites de la famille. Desargues énonce ici une analogie\footnote{p.2, l.10 et suivantes.} entre le cercle et la droite, ce qui peut laisser croire qu'il avait une image juste de ce que l'on appelle aujourd'hui la droite projective (réelle) comme étant topologiquement un cercle et qu'il avait l'idée que les points à l'infini d'un plan formaient en quelque sorte ce que l'on pourrait imaginer être une droite. Le fait que l'ensemble des points à l'infini d'un plan peut être conçu comme une droite est repris plus loin dans le {\it Brouillon} à la fin du développement de sa théorie des \textit{traversales\footnote{Théorie qui deviendra celle, au début du dix-neuvième siècle, des \textit{polaires.}}}, mais nous n'en discuterons pas ici.  


\section{Arbres \& Involutions}
Dans ce qui suit tous les points seront situés sur une même droite~:~il s'agit de géométrie des configurations de points alignés. La première partie concerne la combinatoire pure, qui permet à Desargues de déduire facilement toutes les relations dont il a besoin en n'en démontrant qu'une seule. Ici, l'orientation des segments (ou {\it pièces}) n'a pas d'importance. Dans la deuxième partie nous donnons les définitions d'arbre et d'involution au sens où Desargues l'entend, et démontrons, comme il le fait, que ces deux notions sont équivalentes. Les questions de signe seront ici importantes. 
\subsection{Combinatoire des couples}\label{combinatoire}
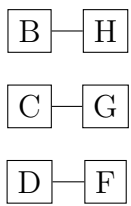
\begin{figure}[!h]
\centering
\begin{tikzpicture}
\node[draw](B) at (0,0) {B};
\node[draw](C) at (0,-1) {C};
\node[draw](D) at (0,-2) {D};
\node[draw](F) at (1,-2) {F};
\node[draw](G) at (1,-1) {G};
\node[draw](H) at (1,0) {H};
\draw (B) -- (H);
\draw (C) -- (G);
\draw (D) -- (F);
\end{tikzpicture}
\caption{Couples de points.}\label{couples}
\end{figure}

Étant donnés trois couples de points distincts $BH, CG, DF$, on dira que $B$ et $H$ ({\it resp.} $C$ et $G$, $D$ et $F$) sont {\it couplés entre eux,} ou simplement {\it couplés\footnote{p.3, l. 20}.} Voici figure \ref{couples} une représentation schématique du couplage choisi par Desargues. Nous la conserverons en mettant sur la même ligne les couples.

\begin{figure}[!h]
\centering
\begin{tikzpicture}
\node[draw](B) at (0,0) {B};
\node[draw](C) at (0,-1) {C};
\node[draw](D) at (0,-2) {D};
\node[draw](F) at (1,-2) {F};
\node[draw,circle](G) at (1,-1) {G};
\node[draw](H) at (1,0) {H};
\draw (G) -- (H);
\draw (G) -- (B);
\draw (G) -- (D);
\draw (G) -- (F);
\end{tikzpicture}
\caption{Les quatre brins de rameau issus de $G$.}\label{brinsderameau}
\end{figure}
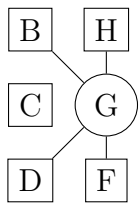

On appelle {\it brin de rameau} toute pièce d'extrémités prises dans deux couples distincts\footnote{p. 3, l. 22}. Le point $G$ par exemple est extrémité de 4 brins de rameau~:~ $GB,GH,GD$ et $GF$. La figure \ref{brinsderameau} en donne une représentation schématique.

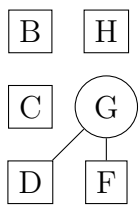
\begin{figure}[!h]
\centering
\begin{tikzpicture}
\node[draw](B) at (0,0) {B};
\node[draw](C) at (0,-1) {C};
\node[draw](D) at (0,-2) {D};
\node[draw](F) at (1,-2) {F};
\node[draw,circle](G) at (1,-1) {G};
\node[draw](H) at (1,0) {H};
\draw (G) -- (D);
\draw (G) -- (F);
\end{tikzpicture}
\caption{La couple de brins de rameau $GD,GF$}\label{couplebrinsderameau}
\end{figure}

Deux brins de rameau forment un {\it couple de brins de rameau} s'ils ont une extrémité commune et si leurs deux autres extrémités sont couplées\footnote{p. 3, l. 27}. En $G$, les brins $GD,GF$ forment un couple\footnote{Nous avons de-ci de-là gardé quelques usages arguésiens, comme l'emploi du féminin pour le mot «~couple~».} de brins de rameau. On dit encore que $GD$ est accouplé à $GF$. La figure \ref{couplebrinsderameau} en donne l'illustration.

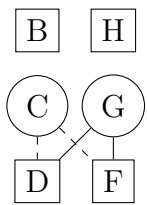
\begin{figure}[!h]
\centering
\begin{tikzpicture}
\node[draw](B) at (0,0) {B};
\node[draw,circle](C) at (0,-1) {C};
\node[draw](D) at (0,-2) {D};
\node[draw](F) at (1,-2) {F};
\node[draw,circle](G) at (1,-1) {G};
\node[draw](H) at (1,0) {H};
\draw (G) -- (D);
\draw (G) -- (F);
\draw[dashed] (C) -- (D);
\draw[dashed] (C) -- (F);
\end{tikzpicture}
\caption{La couple $CD,CF$ est relative de la couple $GD,GF$}\label{couplesrelatifs}
\end{figure}

Étant donné un couple de brins de rameau, on en obtient un autre en changeant l'extrémité commune par le point qui lui est accouplé. Du couple $GD,GF$ on déduit ainsi le couple $CD,CF$. Les deux couples sont dits {\it relatifs\footnote{p. 3, l. 35}.} Ainsi tout couple de brins de rameau admet un unique relatif, comme le montre la figure \ref{couplesrelatifs}.

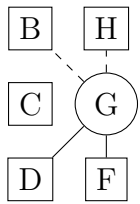
\begin{figure}[!h]
\centering
\begin{tikzpicture}
\node[draw](B) at (0,0) {B};
\node[draw](C) at (0,-1) {C};
\node[draw](D) at (0,-2) {D};
\node[draw](F) at (1,-2) {F};
\node[draw,circle](G) at (1,-1) {G};
\node[draw](H) at (1,0) {H};
\draw[dashed] (G) -- (H);
\draw[dashed] (G) -- (B);
\draw (G) -- (D);
\draw (G) -- (F);
\end{tikzpicture}
\caption{Les deux couples de brins gémeaux issus de $G$.}\label{couplesgemeaux}
\end{figure}

Étant donné un couple de brins de rameau, on en obtient un autre en changeant les extrémités qui sont différentes par les deux extrémités d'un autre couple. Du couple $GD,GF$ on déduit ainsi le couple $GB,GH$. Les deux couples sont dits {\it gémeaux\footnote{p. 3, l. 42}.} Ainsi tout couple de brins de rameau admet un unique gémeau, comme le montre la figure \ref{couplesgemeaux}.

Un couple de brins comme $GD,GF$ donne un rectangle dont les côtés sont $GD, GF$. De la sorte, les notions de relativité et de gémellité s'appliquent aussi aux rectangles définis par les couples.

Desargues définit directement la notion d'arbre ({\it voir} plus loin la définition \ref{definition-arbre}) en ne séparant pas les différents aspects qui sont à l'{\oe}uvre dans sa définition (aspects combinatoires, d'ordonnancement, métriques ou autre). Dans un souci de clarification, nous avons préféré mettre en exergue ces différents aspects en introduisant quelques définitions intermédiaires. 

Si l'on se donne un tiers point $A$, nous disons de la configuration $(A,BH, CG, DF)$ que c'est un {\it proto-arbre} dont $A$ est la {\it souche}, dont les pièces $AB,AH,AC\ldots$ sont les {\it branches,} et dont les points $B,C,D,F,G,H$ sont les {\it n{\oe}uds.} 

\begin{figure}[!h]
\centering
\begin{tikzpicture}
\node[draw,circle](A) at (0.5,1) {A};
\node[draw](B) at (0,0) {B};
\node[draw](C) at (0,-1) {C};
\node[draw](D) at (0,-2) {D};
\node[draw](F) at (1,-2) {F};
\node[draw](G) at (1,-1) {G};
\node[draw](H) at (1,0) {H};
\draw (A) -- (B);
\draw (A) -- (H);
\end{tikzpicture}
\caption{Une couple de branches vers une couple de n{\oe}uds issues de la souche $A$.}\label{couplesbranches}
\end{figure}
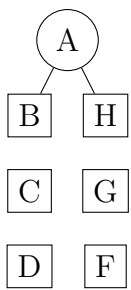

Deux points couplés donnent naissance à un {\it couple de branches\footnote{p. 3, l. 14},} comme le montre la figure \ref{couplesbranches}. Notons tout de suite qu'un brin de rameau est somme ou différence de branches, ce qui sera extensivement exploité par Desargues dans son utilisation de la proposition 12 du livre 5 des \textit{Éléments}, \textit{voir} l'appendice. 
\subsection{Arbre, involution}
Nous nous plaçons maintenant plus directement dans le cadre auquel s'intéresse Desargues, où la droite est {\it orientée,} ou ses points {\it ordonnés.} Donnons-nous un point $A$ et deux segments $AF,AD$ ayant $A$ comme extrémité commune. On dit que $A$ est {\it engagé} dans le couple de branches $AD,AF$ si $A$ est entre les deux points $F$ et $D$\footnote{p. 2, l. 33}.

\begin{figure}[!h]
\centering
\begin{tikzpicture}
\draw (0,0) -- (6,0);
\draw (1,0) node {$\times$};
\draw (1,0) node[above] {$F$};
\draw (4,0) node {$\bullet$};
\draw (4,0) node[above]{$A$};
\draw (5,0) node {$\times$};
\draw (5,0) node[above]{$D$};
\end{tikzpicture}
\caption{La souche $A$ est engagée entre les branches de la couple $AF,AD$.}\label{engage}
\end{figure}
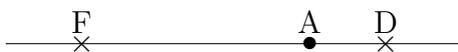

\begin{figure}[!h]
\centering
\begin{tikzpicture}
\draw (0,0) -- (6,0);
\draw (1,0) node {$\bullet$};
\draw (1,0) node[above] {$A$};
\draw (4,0) node {$\times$};
\draw (4,0) node[above]{$D$};
\draw (5,0) node {$\times$};
\draw (5,0) node[above]{$F$};
\end{tikzpicture}
\caption{La souche $A$ est dégagée d'entre les branches de la couple $AF,AD$.}\label{degage}
\end{figure}
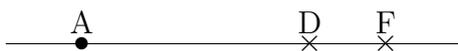

On dit qu'il est {\it dégagé\footnote{p. 2, l. 35}} sinon ({\it voir} les figures \ref{engage} et \ref{degage}).

\begin{definition} Nous appelons \emph{arbre combinatoire} tout proto-arbre $(A,BH,CG,DF)$ tel que $A$ est semblablement engagé ou dégagé à chacun des trois couples de branches $AB,AH;\; AC,AG;\; AD,AF$.
\end{definition}

Ayant bien distingué les aspects combinatoires et d'ordonnancement qui sont au c{\oe}ur de la notion d'arbre\footnote{Disinction que Desargue ne fait pas.}, nous pouvons maintenant donner la première définition importante du {\it Brouillon Project}\footnote{p. 3, l. 5}~:~
\begin{definition}\label{definition-arbre} On appelle \emph{arbre} tout arbre combinatoire $(A,BH,CG,DF)$ qui satisfait à la condition métrique suivante~:~
\[
AB.AH=AC.AG=AD.AF.
\]
\end{definition}

\begin{figure}[!h]
\centering
\begin{tikzpicture}
\node[draw](A) at (0.5,1) {A};
\node[draw](B) at (0,0) {B};
\node[draw](C) at (0,-1) {C};
\node[draw](D) at (0,-2) {D};
\node[draw](F) at (1,-2) {F};
\node[draw](G) at (1,-1) {G};
\node[draw](H) at (1,0) {H};
\node[draw](E) at (0.5,-3) {E$=\infty$};
\draw (A) -- (E);
\draw (B) -- (H);
\draw (C) -- (G);
\draw (D) -- (F);
\end{tikzpicture}
\caption{La souche $A$ comme n{\oe}ud couplé au point à l'infini.}\label{souche-infini}
\end{figure}
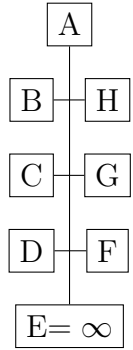

La première chose qui frappe en lisant Desargues (dans les pages 2 à 4 du Brouillon) c'est son choix très soigneux des dénominations des points. Il n'aura pas échappé au lecteur que dans la suite des couples\footnote{Il n'est pas impossible que ce couplage fasse écho au travail de Desargues pour une méthode facile d'apprentissage de la lecture musicale, publié dans l'Harmonie Universelle de Mersenne ({\it voir} \cite{HarmonieUniverselle}, livre sixiesme).} on obtient comme un hexagone dont on nomme les trois premiers sommets dans l'ordre alphabétique $B,C,D$ et les trois suivants dans l'ordre inverse $H,G,F$. Si l'on fait entrer la souche dans le jeu, la lettre $E$ apparaît alors comme la grande absente. Notre thèse est que Desargues a en tête le fait que la souche $A$ apparaîtra comme n{\oe}ud d'un arbre dont l'accouplé est le point à l'infini\footnote{p. 6, l. 33}. On pourrait donc compléter les schémas ci-dessus de la manière illustrée dans la figure \ref{souche-infini}.

Avant d'arriver à cette image, Desargues va montrer qu'une configuration d'arbre donne naissance à une configuration d'où la souche disparaît, ne laissant que les trois couples de points $BH, CG, DF$. Comme pour la notion d'arbre, on peut décomposer la nouvelle notion, dite {\it d'involution,} en la conjonction de deux conditions~:~l'une, combinatoire, traitant de l'ordonnancement des points et l'autre, métrique, traitant de rapports de grandeurs. La combinatoire de l'involution doit prendre en compte la manière dont deux couples de points, ou deux segments, peuvent s'organiser sur une droite ordonnée. Deux couples $CG, DF$ sont dits \emph{mêlés\footnote{p. 2, l. 43}} si chacune des deux pièces $CG, DF$ possède une extrémité dans l'autre et une extrémité hors de l'autre. Deux couples CG, DF sont dits \emph{démêlés\footnote{p. 2, l. 47}} si ou bien chacune des deux pièces a ses extrémités hors de l'autre pièce ou si les deux extémités de l'une des deux pièces sont toutes les deux dans l'autre pièce ({\it voir} les figures \ref{meslez} et \ref{demeslez}).

\begin{figure}[!h]
\centering
\begin{tikzpicture}
\draw (0,0) -- (7,0);
\draw (1,0) node {$\times$};
\draw (1,0) node[above] {$D$};
\draw (3,0) node {$+$};
\draw (3,0) node[above] {$C$};
\draw (4,0) node {$\times$};
\draw (4,0) node[above] {$F$};
\draw (6,0) node {$+$};
\draw (6,0) node[above] {$G$};
\end{tikzpicture}
\caption{Les deux couples $CG,DF$ sont meslez.}\label{meslez}
\end{figure}
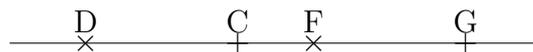
\begin{figure}[!h]
\centering
\begin{tikzpicture}
\draw (0,0) -- (7,0);
\draw (1,0) node {$\times$};
\draw (1,0) node[above] {$D$};
\draw (4,0) node {$+$};
\draw (4,0) node[above] {$C$};
\draw (3,0) node {$\times$};
\draw (3,0) node[above] {$F$};
\draw (6,0) node {$+$};
\draw (6,0) node[above] {$G$};
\end{tikzpicture}

\begin{tikzpicture}
\draw (0,0) -- (7,0);
\draw (3,0) node {$\times$};
\draw (3,0) node[above] {$D$};
\draw (1,0) node {$+$};
\draw (1,0) node[above] {$C$};
\draw (4,0) node {$\times$};
\draw (4,0) node[above] {$F$};
\draw (6,0) node {$+$};
\draw (6,0) node[above] {$G$};
\end{tikzpicture}
\caption{Les deux couples $CG,DF$ sont démeslez.}\label{demeslez}
\end{figure}
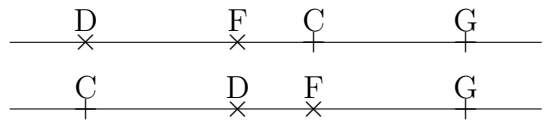

Toujours dans un souci de clarification, nous introduisons la définition suivante~:~
\begin{definition} Nous appelons \emph{involution combinatoire} toute configuration de trois couples de points $(BH, CG, DF)$ telle que soit deux couples quelconques de points sont toujours mêlés, soit deux couples quelconques de points sont toujours démêlés. 
\end{definition}
Cette condition d'ordonnancement, contrairement à celle d'engagement définie pour les arbres, est une notion \emph{projective} qui a un sens sur la droite projective réelle que l'on peut identifier à un cercle. En y rajoutant une condition métrique, on obtient la définition suivante, la plus importante sans doute de tout le {\it Brouillon Project~:}
\begin{definition} On appelle \emph{involution\footnote{p. 4, l. 35}} toute involution combinatoire $(BH, CG, DF)$ satisfaisant à la condition métrique suivante~:~deux rectangles relatifs sont entre eux comme leurs gémeaux. Dit autrement~:~
\[
\frac{GF.GD}{CF.CD}=\frac{GB.GH}{CB.CH},\; \frac{FC.FG}{DC.DG}=\frac{FB.FH}{DB.DH},\; \frac{HC.HG}{BC.BG}=\frac{HD.HF}{BD.BF}.
\]
\end{definition}
La figure \ref{involution} permet de mieux comprendre l'interaction entre la combinatoire et les rapports métriques.

\begin{figure}[!h]
\centering
\begin{tikzpicture}
\draw (1,0) node[above] {$GF.GD$};
\draw (0,0) -- (2,0);
\draw (1,0) node[below] {$CF.CD$};
\draw [<->] (0,-0.3) to[bend left] (0,0.3);
\draw (-1,0) node {relatifs};
\draw (2.2,0) node {$=$};
\draw (2.5,0) -- (4.5,0);
\draw (3.5,0) node[above]{$GB.GH$};
\draw (3.5,0) node[below]{$CB.CH$};
\draw [<->] (4.5,-0.3) to[bend right] (4.5, 0.3);
\draw (5.4,0) node {relatifs};
\draw  [<->] (1,0.7) to[bend left] (3.5,0.7);
\draw (2.3, 1.4) node {gémeaux};
\end{tikzpicture}
\caption{La combinatoire~:~rapports de rectangles relatifs égaux aux rapports de leurs gémeaux.}\label{involution}
\end{figure}

La figure \ref{involution-m} reprend ces mêmes rapports avec un emploi plus moderne des lettres, de manière à permettre au lecteur d'assimiler plus facilement la combinatoire gemellaire.

\begin{figure}[!h]
\centering
\begin{tikzpicture}
\draw (1,0) node[above] {$XY.XY'$};
\draw (0,0) -- (2,0);
\draw (1,0) node[below] {$X'Y.X'Y'$};
\draw [<->] (0,-0.3) to[bend left] (0,0.3);
\draw (-1,0) node {relatifs};
\draw (2.2,0) node {$=$};
\draw (2.5,0) -- (4.5,0);
\draw (3.5,0) node[above]{$XZ.XZ'$};
\draw (3.5,0) node[below]{$X'Z.X'Z'$};
\draw [<->] (4.5,-0.3) to[bend right] (4.5, 0.3);
\draw (5.4,0) node {relatifs};
\draw  [<->] (1,0.7) to[bend left] (3.5,0.7);
\draw (2.3, 1.4) node {gémeaux};
\end{tikzpicture}
\caption{La combinatoire~:~rapports de rectangles relatifs égaux aux rapports de leurs gémeaux pour l'involution (XX',YY',ZZ').}\label{involution-m}
\end{figure}

\subsection{L'équivalence arbre-involution}
Desargues ne procède bien évidemment pas comme nous venons de le faire~:~il définit les termes de tronc, de n{\oe}ud, de rameau, de brin ainsi que ceux d'engagement page 2, lignes 16 à 35 du {\it Brouillon.} La notion d'emmêlement est définie sur la même page, un peu plus bas, lignes 43 à 46. Il donne alors la définition d'un arbre, en combinant les conditions d'ordonnancement et de métrique, page 3, ligne 5. Après avoir défini les termes arboricoles, il déploie la combinatoire décrite dans la section \ref{combinatoire} aux lignes 20 à 46 de la page 3. Il commence alors à démontrer que d'un arbre, on tire une involution, terme qu'il ne définira qu'à la fin de sa démonstration. Énonçons le théorème démontré par Desargues~:~
\begin{theoreme}\label{arbre-involution} Si $(A;BH, CG, DF)$ est un arbre, alors $(BH, CG, DF)$ est une involution.
\end{theoreme}
Là encore il ne donne pas explicitement un tel énoncé (puisqu'en particulier il n'a pas défini la notion d'involution) mais se lance directement, page 3, ligne 47, dans une série de considérations qui en constitue une preuve. Notons enfin que, de manière implicite, Desargues ne traite que le cas de configurations «~génériques~», quand par exemple tous les points sont distincts deux à deux, de manière à ce que toutes ses manipulations soient licites. Les cas particuliers, plus précisément ceux où des branches couplées sont égales, seront traités plus loin par Desargues dans son long développement sur les n{\oe}uds moyens.

La preuve du théorème \ref{arbre-involution} par Desargues se fait en termes exclusivement euclidiens et c'est celle, peu ou prou, que nous donnons maintenant.
\begin{demonstration} Desargues commence par passer rapidement sur le problème de l'agencement combinatoire des points. Une difficulté qu'il rencontrera dans la preuve de la réciproque du théorème  est que la condition d'emmêlement n'est conséquence de la condition d'engagement de la souche {\it que sous l'hypothèse que la condition métrique est aussi vérifiée.} Plus exactement, il souligne que si la condition métrique sur les rectangles d'un arbre a lieu, alors les conditions d'engagement de la souche et d'emmêlement des couples sont équivalentes. Cela est fait entre la ligne 47 de la page 3 et la ligne 2 de la page 4, et est clairement énoncé aux lignes 3 à 5 de la page 4. Donnons cette démonstration et supposons donc que $AB.AH=AC.AG=AD.AF$. 
\begin{itemize}
\item Si l'on est dans la situation où la souche est engagée on peut supposer, sans perte de généralité, que $F$ et $C$ sont du même côté de la souche et que $F$ est plus éloigné de la souche que ne l'est $C$. Mais alors la condition sur les rectangles impose que $D$ soit plus près de la souche que $G$ ne l'est et donc les deux couples $CG, DF$ sont mêlés. On voit donc que si la souche est semblablement engagée, alors les couples sont semblablement mêlés, {\it voir} la figure \ref{engages-meslez}.

\begin{figure}[!h]
\centering
\begin{tikzpicture}
\draw (0,0) -- (7,0);
\draw (1,0) node {$\times$};
\draw (1,0) node[above] {$F$};
\draw (3,0) node {$+$};
\draw (3,0) node[above] {$C$};
\draw (3.5,0) node[above] {$A$};
\draw (3.5,0) node {$\bullet$};
\draw (4,0) node {$\times$};
\draw (4,0) node[above] {$D$};
\draw (6,0) node {$+$};
\draw (6,0) node[above] {$G$};
\end{tikzpicture}
\caption{De souche engagée à couples meslez.}\label{engages-meslez}
\end{figure}
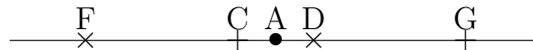

\item Si l'on est dans la situation où la souche est dégagée on peut supposer, sans perte de généralité, que $F$ est plus proche de la souche que ne le sont $C,D$ et $G$. Mais alors la condition sur les rectangles impose que $D$ soit plus éloigné de $A$ que $C$ et $G$ ne le sont et donc les deux couples $CG, DF$ sont démêlés. On voit donc que si la souche est semblablement dégagée, alors les couples sont semblablement démêlés, {\it voir} la figure \ref{degages-demeslez}\footnote{Remarquons que si $C$ est de l'autre côté de la souche par rapport à $F$, alors $G$ l'est aussi et l'on est encore en situation démêlée.}.

\begin{figure}[!h]
\centering
\begin{tikzpicture}
\draw (0,0) -- (7,0);
\draw (1,0) node {$\bullet$};
\draw (1,0) node[above] {$A$};
\draw (2,0) node {$\times$};
\draw (2,0) node[above] {$F$};
\draw (3.5,0) node[above] {$C$};
\draw (3.5,0) node {$+$};
\draw (4,0) node {$+$};
\draw (4,0) node[above] {$G$};
\draw (6,0) node {$\times$};
\draw (6,0) node[above] {$D$};
\end{tikzpicture}
\caption{De souche dégagée à couples démeslez.}\label{degages-demeslez}
\end{figure}
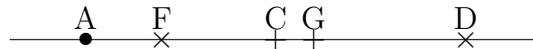

\item Réciproquement, si les deux couples $CG, DF$ sont mêlés, alors la souche ne peut être dégagée des deux couples. En effet, si cela était le cas, nous pouvons supposer, sans perte de généralité, que $F$ est plus proche de $A$ que $C$ ne l'est. Mais comme $CG$ et $DF$ sont mêlés, $D$ est entre $C$ et $G$ et donc $AF.AD<AC.AG$, contredisant la condition métrique sur les rectangles. On voit donc que si les couples sont semblablement mêlés, alors la souche est semblablement engagée.
\item Si les deux couples $CG, DF$ sont démêlés, alors la souche ne peut être engagée. Si cela était le cas, nous pourrions supposer, sans perte de généralité, que $C$ et $F$ sont du même côté de la souche et que $C$ est plus loin de $A$ que $F$ ne l'est. Les points $D$ et $G$ sont alors de l'autre côté de $A$ et comme les couples sont démêlés, $G$ est plus loin de $A$ que $D$ ne l'est, entrainant que $AF.AD<AC.AG$ et contredisant la condition métrique sur les rectangles. On voit donc que si les couples sont semblablement démêlés, alors la souche est semblablement dégagée.
\end{itemize}
Revenons maintenant à la démonstration du théorème proprement dit. Comme la configuration est un arbre, la condition combinatoire de semblable emmêlement est vérifiée et il nous reste uniquement à démontrer l'égalité des rapports de rectangles relatifs. Cette partie de la démonstration occupe Desargues des lignes 7 à 34 de la page 4 du {\it Brouillon Project.} Prenons deux couples quelconques, par exemples $CG,DF$. De l'égalité de rectangles $AC.AG=AD.AF$ on tire les deux proportions\footnote{Nous employons le langage euclidien~:~quatre grandeurs $a,b,c,d$ sont proportionnelles si $a/b=c/d$. Un telle égalité de rapports s'appelle une \textit{proportion.}}
\[
\frac{AC}{AD}=\frac{AF}{AG}\;\mbox{et}\; \frac{AC}{AF}=\frac{AD}{AG}.
\]
Si nous sommes dans la situation de la figure \ref{engages-meslez}, on peut écrire l'égalité de longueurs\footnote{Nous ne considérons ici que des longueurs ou grandeurs positives, comme le fait Desargues. Nous renvoyons le lecteur à la section suivante pour une version algébrique plus moderne de cette preuve.} $AF-AC=CF$ et $AG-AD=GD$ et la première proportion permet alors d'écrire que
\[
\frac{AC}{AD}=\frac{AF}{AG}=\frac{AF-AC}{AG-AD}=\frac{CF}{GD}.
\]
Si nous sommes dans la situation de la figure \ref{degages-demeslez}, on peut écrire l'égalité de longueurs $AC-AF=CF$ et $AD-AG=GD$ et la première proportion permet alors d'écrire que
\[
\frac{AC}{AD}=\frac{AF}{AG}=\frac{AC-AF}{AD-AG}=\frac{CF}{GD},
\]
résultat semblable au précédent. 

Si enfin nous étions dans la situation démêlée où $D$ et $F$ sont du côté opposé de $C$ et $G$ par rapport à $A$, nous tirerions de la première proportion la même identité de rapport~:~
\[
\frac{AC}{AD}=\frac{AC+AF}{AD+AG}=\frac{CF}{GD}.
\]
De la seconde proportion on déduit de la même manière~:~
\[
\frac{AC}{AF}=\frac{AD}{AG}=\frac{CD}{GF}.
\]
Comme l'on veut faire apparaître des rapports de rectangles relatifs, on écrit le rapport $AC/AG$ comme une composition de rapports habilement choisis~:~
\[
\frac{AC}{AG}=\frac{AC}{AD}\frac{AD}{AG}=\frac{CF}{GD}\frac{CD}{GF}=\frac{CF.CD}{GD.GF}.
\]
Remarquons que l'autre composition de rapports
\[
\frac{AC}{AG}=\frac{AC}{AF}\frac{AF}{AG}
\]
donnerait le même résultat. 

Nous avons donc démontré la formule de la ligne 15 de la page 4, à savoir~:~
\begin{equation}\label{ACAG}
\frac{AC}{AG}=\frac{CD.CF}{GD.GF},
\end{equation}
où apparaît un rapport de rectangles relatifs (c-à-d. dont les côtés sont des couples de brins relatifs entre eux). On peut associer à l'égalité \ref{ACAG} le diagramme  de la figure \ref{ACAG-DF}.

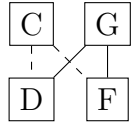
\begin{figure}[!h]
\centering
\begin{tikzpicture}
\node[draw](C) at (0,-1) {C};
\node[draw](D) at (0,-2) {D};
\node[draw](F) at (1,-2) {F};
\node[draw](G) at (1,-1) {G};
\draw (G) -- (D);
\draw (G) -- (F);
\draw[dashed] (C) -- (D);
\draw[dashed] (C) -- (F);
\end{tikzpicture}
\caption{Le diagramme de l'identité $AC/AG=(CD.CF)/(GD.GF)$}\label{ACAG-DF}
\end{figure}

Pour faire entrer en jeu les rectangles gémeaux, gardons le couple $CG$ et prenons pour autre couple $BH$. De l'identité de rectangles $AB.AH=AC.AG$ on tire les deux proportions
\[
\frac{AC}{AB}=\frac{AH}{AG}\;\mbox{et}\; \frac{AC}{AH}=\frac{AB}{AG}.
\]
En utilisant les mêmes arguments que ci-dessus, on démontre les deux égalités
\[
\frac{AC}{AB}=\frac{AH}{AG}=\frac{CH}{GB}\;\mbox{et}\; \frac{AC}{AH}=\frac{AB}{AG}=\frac{CB}{GH}.
\]
On peut alors reprendre le rapport $AC/AG$ et l'écrire comme composé en utilisant cette fois-ci un des points du nouveau couple $BH$, par exemple $H$~:~
\[
\frac{AC}{AG}=\frac{AC}{AH}\frac{AH}{AG}=\frac{CB.CH}{GB.GH},
\]
identité dont le diagramme associé est donné dans la figure \ref{ACAG-BH} et d'où découle, avec l'identité démontrée plus haut, que
\[
\frac{CD.CF}{GD.GF}=\frac{CB.CH}{GB.GH}.
\]

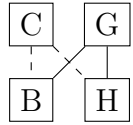
\begin{figure}[!h]
\centering
\begin{tikzpicture}
\node[draw](C) at (0,-1) {C};
\node[draw](B) at (0,-2) {B};
\node[draw](H) at (1,-2) {H};
\node[draw](G) at (1,-1) {G};
\draw (G) -- (B);
\draw (G) -- (H);
\draw[dashed] (C) -- (B);
\draw[dashed] (C) -- (H);
\end{tikzpicture}
\caption{Le diagramme de l'identité $AC/AG=(CB.CH)/(GB.GH)$}\label{ACAG-BH}
\end{figure}

Desargues répète plusieurs fois les arguments détaillés ci-dessus, aux lignes 15 à 28 de la page 4. Toujours est-il que l'on démontre ainsi que deux branches couplées comme $AC$ et $AG$ sont dans le même rapport que deux rectangles relatifs quelconques basés sur ces deux branches, rapport qui est bien sûr préservé par passage aux gémeaux. On peut donc faire tourner aveuglément la combinatoire et ainsi obtenir les autres égalités de rapports. Les couples jouant des rôles symétriques, Desargues pourrait s'arrêter là et invoquer la généralité de l'argument qui précède. Mais il insiste sur l'utilisation de sa combinatoire~:~ce qui apparaît de prime abord comme un simple outil technique permettant de faciliter la permutation des lettres est en fait l'une des véritables nouveautés de la démarche arguésienne. Ainsi page 4, ligne 29, il obtient
\[
\frac{AF}{AD}=\frac{FC.FG}{DC.DG}=\frac{FB.FH}{DB.DH},
\]
dont les diagrammes sont donnés par la figure \ref{AFAD}.

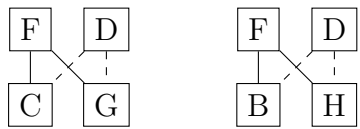
\begin{figure}[!h]
\centering
\begin{tikzpicture}
\node[draw](F) at (0,0) {F};
\node[draw](C) at (0,-1) {C};
\node[draw](D) at (1,0) {D};
\node[draw](G) at (1,-1) {G};
\draw (F) -- (C);
\draw (F) -- (G);
\draw[dashed] (D) -- (C);
\draw[dashed] (D) -- (G);
\node[draw](F1) at (3,0) {F};
\node[draw](C1) at (3,-1) {B};
\node[draw](D1) at (4,0) {D};
\node[draw](G1) at (4,-1) {H};
\draw (F1) -- (C1);
\draw (F1) -- (G1);
\draw[dashed] (D1) -- (C1);
\draw[dashed] (D1) -- (G1);
\end{tikzpicture}
\caption{Diagrammes des identités $AF/AD=(FC.FG)/(DC.DG)=(FB.FH)/(DB.DH)$}\label{AFAD}
\end{figure}

De même à la ligne 32 de la page 4 il obtient
\[
\frac{AH}{AB}=\frac{HC.HG}{BC.BG}=\frac{HD.HF}{BD.BF},
\]
dont les diagrammes sont donnés par la figure \ref{AHAB}, achevant ainsi la démonstration du théorème. 
\end{demonstration}
\begin{figure}[!h]
\centering
\begin{tikzpicture}
\node[draw](H) at (0,0) {H};
\node[draw](C) at (0,-1) {C};
\node[draw](B) at (1,0) {B};
\node[draw](G) at (1,-1) {G};
\draw (H) -- (C);
\draw (H) -- (G);
\draw[dashed] (B) -- (C);
\draw[dashed] (B) -- (G);
\node[draw](H1) at (3,0) {H};
\node[draw](C1) at (3,-1) {D};
\node[draw](B1) at (4,0) {B};
\node[draw](G1) at (4,-1) {F};
\draw (H1) -- (C1);
\draw (H1) -- (G1);
\draw[dashed] (B1) -- (C1);
\draw[dashed] (B1) -- (G1);
\end{tikzpicture}
\caption{Diagrammes des identités $AH/AB=(HC.HG)/(BC.BG)=(HD.HF)/(BD.BF)$}\label{AHAB}
\end{figure}
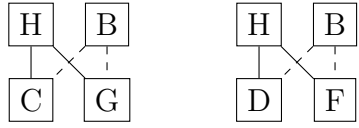

Le lecteur attentif remarquera à nouveau que les arguments avancés ci-dessus ne sont valables que si d'une part aucun n{\oe}ud de l'arbre n'est confondu avec la souche (si $B=A$ on ne peut diviser par $AB$) et si d'autre part deux n{\oe}uds de deux couples différents sont eux aussi différents (si $D=C$ on ne peut diviser par $DC.DG$). Ces cas particuliers seront traités un peu plus tard par Desargues. Celui d'un n{\oe}ud confondu avec la souche permettra de montrer que la souche elle-même peut être vue comme accouplée au point à l'infini, et celui de n{\oe}uds confondus sera l'occasion de développer la théorie des {\it n{\oe}uds moyens simples et doubles.} 
\subsubsection{D'une involution à un arbre~:~la situation générique}
Après avoir démontré le théorème \ref{arbre-involution}, Desargues donne la définition de la notion d'involution, lignes 35 à 40 du {\it Brouillon Project.} La suite du texte est difficile à mettre en ordre. Aux lignes 46 et suivantes de la page 4, il répète que la notion d'involution sous-entend celle d'arbre. À partir de la ligne 47 suit un argument qui ressemble à une contraposition ou à un raisonnement par l'absurde, sur lequel nous reviendrons. Les choses paraissent se préciser ligne 12 de la page 5, où Desargues semble affirmer que si $GC,DF$ est couple de n{\oe}uds, alors la souche $A$ est donnée de position. Il s'embarque alors dans une série d'arguments (lignes 16 à 33), qu'il répète plusieurs fois, tentant de montrer la {\it réciproque} du théorème \ref{arbre-involution}. Il n'y parvient pas réellement mais on peut en reconstruire une preuve à partir des éléments désordonnés qu'il présente.
\begin{theoreme}\label{involution-arbre} Si $(BH,CG,DF)$ est une involution, il existe un point $A$ tel que $(A;BH,CG,DF)$ soit un arbre.
\end{theoreme}
Les raisons de la confusion de la présentation de Desargues tiennent à la difficulté d'organiser clairement la démonstration lorsque l'on n'a pas séparé les diverses conditions qui définissent arbres et involutions. De ce fait envisager une preuve par l'absurde ou par contraposée semble difficile. Les arguments de Desargues semblent s'enchaîner sans ordre entre la fin de la page 4 et la fin de la page 5. Nous verrons que dans sa lettre, Beaugrand présente une preuve très lisible du théorème \ref{arbre-involution}, mais que, comme Desargues, il prétend que sa réciproque est évidente en faisant appel au raisonnement par contraposée. En ordonnant quelque peu les arguments de Desargues, on peut cependant reconstruire une preuve de cette réciproque. 
\begin{demonstration} Supposons que $(BH,CG,DF)$ est une involution et supposons de plus que les six points sont distincts deux-à-deux de manière à pouvoir considérer tous les rapports dont nous aurons besoin. Concentrons-nous sur deux couples arbitraires $CG,DF$. On peut construire un point $A$ et un seul\footnote{L'unicité se déduit de la condition d'engagement.} tel que le rapport de $AG$ à $AF$ soit le même que celui de $GD$ à $FC$. En termes algébriques modernes, il s'agit de résoudre une équation du premier degré, ce que nous verrons plus loin. Le point $A$ sera alors semblablement engagé ou dégagé aux deux branches de chaque couple suivant que les couples sont mêlés ou démêlés, comme nous l'avons vu au début de la preuve du théorème \ref{arbre-involution}. 

Le même argument euclidien sur les proportions que celui utilisé dans la preuve du théorème \ref{arbre-involution} nous permet de dire que
\begin{equation}\label{prout1}
\frac{AG}{AF}=\frac{GD}{FC}=\frac{AD}{AC}
\end{equation}
d'où\footnote{p.5, l. 21--23} $AG.AC=AD.AF$. L'on en déduit alors\footnote{p.5, l. 24} que
\begin{equation}\label{prout2}
\frac{CF.CD}{GD.GF}=\frac{AC}{AG}.
\end{equation}
Ensuite Desargues, aux lignes 27 à 33, se répète, mais ne démontre pas ce qui manque encore pour prouver que $A$ est souche de l'arbre $(A;BH,CF,DG)$, à savoir que $AB.AH=AC.AG=AD.AF$. Or cela peut se déduire d'un argument qu'il a employé plus haut dans le {\it Brouillon,} de la ligne 47 de la page 4 à la ligne 11 de la page 5. Il s'agit de raisonner par l'absurde en supposant que $AB.AH\neq AC.AG$. Comme dit plus haut, on peut construire un point $H'$ tel que $AB.AH'=AC.AG$ et donc tel aussi que $AB.AH'=AD.AF$ d'après \ref{prout1}. De ces égalités de rectangles on tire, comme on l'a fait plus haut, des égalités de rapports~:~
\[
\frac{AC}{AG}=\frac{CB.CH'}{GB.GH'}.
\]
Du fait que $H\neq H'$ on déduit, grâce au théorème de Thalès et au cinquième postulat d'Euclide, que
\[
\frac{CH'}{GH'}\neq \frac{CH}{GH}
\]
ce qui, combiné avec l'identité \ref{prout2} permet d'écrire que
\[
\frac{CF.CD}{GD.GF}\neq \frac{CB.CH}{GB.GH},
\]
contredisant le fait que $(BH,CG,DF)$ est une involution et achevant ainsi la démonstration du théorème \ref{involution-arbre}.
\end{demonstration}
\subsubsection{Arbres et involutions : la situation dégénérée, n{\oe}uds moyens}
Dans sa définition de la notion d'arbre page 4 et sa démonstration du lien avec la notion d'involution en fin de page 5, Desargues se place implicitement dans une situation \emph{générique,} au sens où les points sont tous distincts deux-à-deux et où toutes les opérations qu'il effectue sont licites. À partir de la ligne 41 de la page 5, dans un long développement qui s'étendra jusqu'à la page 10, il s'attaque à l'étude de cas particuliers, dont l'importance tient à leur apparition naturelle dans l'étude des coniques~:~ceux où certains points viennent à se confondre, du fait par exemple de l'égalité de longueur entre brins couplés. La notion d'involution étant de nature projective, le lecteur moderne comprend qu'ici l'auteur risque de buter sur quelques difficultés liées à la confusion des structures.  

De prime abord, le texte de Desargues semble confus. Cela tient à son choix, dont il a averti le lecteur en ouverture, de ne pas faire «~la distinction necessaire d'entre les impositions de nom, autrement definitions, les propositions, les demonstrations, quand elles sont en suitte\footnote{p.~1, l.~1 à 4}.~» Il commence par donner les définitions de n{\oe}ud moyen, simple et double, ainsi que de n{\oe}ud extrême, intérieur ou extérieur. Il passe ensuite à l'examen de ce qui se passe lorsque l'un des n{\oe}uds vient à se confondre avec la souche, passage qui est indépendant de la notion de n{\oe}ud moyen, mais qui s'avèrera cependant important pour l'étude de ceux-ci. Vient ensuite une étude du cas des n{\oe}uds moyens simples, situation qu'il évacue en ne la considérant pas encore comme une involution~:~«~Il y a nombre d'autres proprietez particulieres à ce cas de cette espece de conformation d'arbre, où chacun peut s'égayer à sa fantaisie, mais il n'est pas encore icy du nombre de ceux qui constituent une involution\footnote{p.~7, l.~23}.~» Il s'attaque alors aux n{\oe}uds moyens doubles, dont il développe l'étude plus longuement. Il y montrera notamment qu'apparaît dans cette situation une division harmonique, dont il montrera qu'elle peut être considérée de deux manières comme une involution, en choisissant convenablement des souches pour les arbres correpondants~:~c'est sa notion de \emph{souches réciproques entre elles.} Ayant clarifié la combinatoire, il se lance dans la démonstration d'un grand nombre de relations de proportions, que nous détaillons plus bas et qui ont dû rebuter bien des lecteurs.

La synthèse se trouve finalement dans un paragraphe sybillin de la page 10, aux lignes 31 à 36~:~«~Quand en un tronc droict, trois couples de n{\oe}uds extremes DF, CG, BH, sont en involution entre eux, que deux autres couples de n{\oe}uds, moyens, unis, doubles ou simples, PQ,XY, sont une involution de quatre points avec chacune des deux quelconques couples CG,\& BH, de ces trois couples de n{\oe}uds extremes. Ces deux mesmes n{\oe}uds moyens PQ, XY, sont encore une involution aussi de quatre points avec la troisiéme de ces couples de n{\oe}uds extremes DF.~» 

Dit autrement, si $(DF,CG,BH)$ est une involution de six points {\it en position générale,} et si par ailleurs $(PQ,XY,CG)$ ainsi que $(PQ,XY,BH)$ forment des involutions, que $P,Q,X,Y$ soient n{\oe}uds moyens simples ou moyens doubles, peu importe, la configuration $(PQ,XY,DF)$ est une involution (de quatre points seulement, {\it voir} plus bas). C'est à notre sens le premier énoncé d'importance du {\it Brouillon--Project} et à lui seul il justifie l'introduction de cette nouvelle notion qu'est l'involution, qui est \emph{agrégative} au sens où nous l'avons entendu dans l'introduction.

Afin de simplifier la lecture, nous allons traiter ces cas les uns après les autres, en respectant peu ou prou l'ordre désordonné de Desargues.

\vs
\noindent{\bf Les {\it noms imposez}, ou définitions.} Supposons que $(A;BH,CG,DF)$ est un arbre. Une première forme de non-généricité peut être celle où les longueurs des branches d'un même couple sont égales. Si $AC=AG$, on dit que le couple de n{\oe}uds $C,G$ forme un couple de {\it n{\oe}uds moyens\footnote{p. 3, l.12, pour les branches plutôt que pour les n{\oe}uds correspondants.}.} Si la souche $A$ est engagée entre les pièces $AC,AG$, alors $A$ est le milieu de $CG$ de sorte que $C\neq G$. Desargues dit que dans ce cas les n{\oe}uds $C,G$ sont des n{\oe}uds moyens {\it simples\footnote{p.5, l. 41}.} Si au contraire la souche $A$ est dégagée, alors $C=G$ et le couple $C,G$ est un couple de n{\oe}uds moyens {\it doubles\footnote{p. 5, l. 56}.} Enfin, si $AB\neq AH$, on dit que le couple $B,H$ forme un couple de n{\oe}uds {\it extrêmes\footnote{p. 6, l. 8}.} 

Dans la situation où l'arbre ne comporte qu'un couple de n{\oe}uds moyens et deux couples de n{\oe}uds extrêmes, tout se passe bien dans les arguments développés jusque là~:~la notion d'arbre est équivalente à celle d'involution. C'est dans les cas où l'arbre compte {\it deux couples de n{\oe}uds moyens} que la situation se gâte. Faisons ici une remarque que Desargues ne daigne pas écrire~:~un arbre comporte au plus deux couples de n{\oe}uds moyens. Plus précisément même, étant donné un arbre $(A;BH,CG,DF)$ avec $B,H$ n{\oe}uds extrêmes, on peut toujours trouver deux couples et deux seulement $C',G'$ et $D',F'$ de sorte que $(A;BH,C'G',D'F')$ soit un arbre dont $C',G'$ et $D',F'$ soient n{\oe}uds moyens. 

Plaçons-nous dorénavant dans le cas où $(A;BH,CG,DF)$ est un arbre ayant deux couples de n{\oe}uds moyens $C,G$ et $D,F$ et un couple de n{\oe}uds extrêmes $B,H$.

\vs
\noindent{\bf Cas des n{\oe}uds partant à l'infini.} Commençons par examiner ce qui se passe quand l'un des n{\oe}uds extrêmes approche la souche. Que les n{\oe}uds moyens soient simples ou doubles, l'un des n{\oe}uds extrêmes, $B$ par exemple, est entre les n{\oe}uds moyens et l'autre, $H$, ne l'est pas. On dit que $B$ est un n{\oe}ud extrême {\it intérieur}\footnote{p. 6, l. 13} et $H$ un n{\oe}ud extrême {\it extérieur\footnote{p. 6, l. 15}.} De l'égalité de rectangles $AB.AH=AC.AG=AC^2$, il découle qu'à mesure que $B$ approche de la souche $A$, le point $H$ s'en éloigne jusqu'à ce que Desargue appelle la «~distance infinie.~» Cette remarque lui permet, page 6, lignes 33 à 35, de considérer que le point $A$, accouplé au point à distance infinie, que Desargues ne nomme pas mais que nous avons plus haut proposé de noter $E$, forme un couple de points dans l'involution. Il est curieux que ce développement se fasse au milieu de celui concernant les n{\oe}uds simples, puisqu'il en est en fait indépendant.

\vs
\noindent{\bf N{\oe}uds moyens simples.} À partir de la moitié de la page 6 et jusqu'à la page 7, Desargues étudie la situation suivante~:~considérons un arbre $(A;BH, CG, DF)$ dont la souche $A$ est engagée entre deux couples de n{\oe}uds moyens simples de sorte que par exemple $C\neq G$ et $D\neq F$, tandis que $BH$ forme un couple de n{\oe}uds extrêmes\footnote{p. 6, l. 38}. Les identités issues de la définition d'un arbre sont rappellées par Desargues et sont en l'occurence
\[
AB.AH=AC.AG=AD.AF=AC^2=AF^2=AD^2=AG^2.
\]
Il souligne alors\footnote{p. 6, l. 47} que cela signifie en termes euclidiens qu'une quelconque des branches moyennes, $AC$ par exemple, est moyenne proportionnelle des branches extêmes $AB, AH$. On peut sans doute émettre l'hypothèse que c'est cette occurence d'une situation euclidienne classique qui a incité Desargues à utiliser les termes de n{\oe}uds {\it moyens} et {\it extrêmes,} renforçant l'idée que ses choix en ce domaine sont très réfléchis\footnote{Rappelons que si $a/c=c/b$, on dit que $c$ est moyenne proportionnelle des deux extrêmes $a$ et $b$.}.

Comme $A$ est le milieu de $CG$ et de $DF$ qui sont tous deux de longueurs égales, il est nécessaire que, par exemple, $C=F$ et $D=G$, comme il est montré dans la figure 6, deuxième ligne, p.~115 de \cite{taton}. Desargues va démontrer l'égalité suivante~:~
\[
GB.GH=CB.CH
\]
et ce de trois manières différentes. La première démonstration qu'il donne est fautive, mais instructive\footnote{p. 6, l. 49 à l. 56}. Comme les points sont disposés comme en involution, on a
\[
\frac{GD.GF}{CD.CF}=\frac{GB.GH}{CB.CH},
\]
ce qui est étrange puisqu'ici $G=D$ et $C=F$ de sorte que le premier rapport n'est autre que $0/0$. Cela n'effraie pas Desargues qui continue et en déduit immédiatement que 
\[
\frac{GD.GF}{GB.GH}=\frac{CD.CF}{CB.CH},
\]
et comme\footnote{p. 6, l. 54} $GD.GF=CD.CF$, que Desargues déduit probablement de l'égalité $GD=CF$ (égaux en fait à $0$) et de $GF=CD$, on en tire l'identité recherchée, à savoir $GB.GH=CB.CH$.

Il passe ensuite au second argument, annonçant\footnote{p. 6, l. 57} qu'il est évident que
\[
\frac{CH}{BG}=\frac{GH}{BC},
\]
d'où l'identité recherchée suit immédiatement. L'identité ci-dessus est encore une fois conséquence de la définition d'un arbre. De $AB.AH=AC.AG$ on tire
\[
\frac{AH}{AG}=\frac{AC}{AB}
\]
d'où
\[
\frac{AH}{AG}=\frac{AH-AC}{AG-AB}=\frac{AH-AG}{AC-AB}=\frac{GH}{BC}
\]
et
\[
\frac{AH}{AG}=\frac{AH+AC}{AG+AB}=\frac{CH}{BG}.
\]

Son troisième argument est donné à partir de l'avant dernière ligne de la page 6 et se poursuit page 7. Comme il l'a déjà démontré, $AG$ est à son accouplée $AC$ comme le rectangle $GB.GH$ est à son relatif le rectangle $CB.CH$ ({\it voir} la figure \ref{ACAG-BH}). Mais ici $AG=AC$ puisque $CG$ est un couple de n{\oe}uds moyens. Ainsi a-t-on nécessairement $GB.GH=CB.CH$. Et Desargues de conclure ces trois arguments\footnote{p. 7, l2} par la remarque que cette identité devient incompréhensible lorsque $B$ s'unit à la souche $A$ et que le point $H$ part à l'infini. 

Il conclut cette étude par une réinterprétation de la configuration en termes de trois segments $HG, GB, BC$ consécutifs. 

\begin{figure}[!h]
\centering
\definecolor{xdxdff}{rgb}{0.49019607843137253,0.49019607843137253,1.}
\definecolor{qqqqff}{rgb}{0.,0.,1.}
\definecolor{ffxfqq}{rgb}{1.,0.4980392156862745,0.}
\begin{tikzpicture}[line cap=round,line join=round,>=triangle 45,x=1.0cm,y=1.0cm]
\clip(1.26,2.7) rectangle (15.9,5.16);
\draw (2.12,4.18)-- (15.08,4.18);
\begin{scriptsize}
\draw [fill=ffxfqq] (2.12,4.18) circle (2.5pt);
\draw[color=ffxfqq] (2.02,3.79) node {$C$};
\draw [fill=qqqqff] (15.08,4.18) circle (2.5pt);
\draw[color=qqqqff] (15.02,3.87) node {$H$};
\draw [fill=xdxdff] (4.376801221223264,4.18) circle (2.5pt);
\draw[color=xdxdff] (4.38,3.85) node {$B$};
\draw [fill=ffxfqq] (9.868241120372637,4.18) circle (2.5pt);
\draw[color=ffxfqq] (9.82,3.79) node {$G$};
\end{scriptsize}
\end{tikzpicture}
\end{figure}

De l'égalité de rectangles $GB.GH=CB.CH$ il déduit que
\[
\frac{CH}{GB}=\frac{GH}{CB},
\]
ou encore que la somme des trois segments est au segment du milieu comme le segment du bout, du côté du n{\oe}ud extrême extérieur $H$, est au segment de l'autre bout~:~
\[
\frac{CB+BG+GH}{BG}=\frac{GH}{CB}.
\]
Le cas où le n{\oe}ud extrême extérieur est à l'infini\footnote{p. 7, l. 20}, c'est-à-dire où $B=A$, donne un cas particulier intéressant que l'on peut déduire de la configuration d'arbre mais aussi, directement, des considérations sans souche faites juste avant~:~$B$ devient milieu de $CG$. Desargues finit par dire\footnote{p. 7, l. 23} qu'en cette espèce de conformation il y a nombre d'autres propriétés intéressantes, mais qu'il ne la considère pas ici encore\footnote{l. 24} comme constituant une involution, ce qui est quelque peu surprenant mais est éclairci plus loin, page 10, lignes 31 et suivantes, dans le paragraphe sybillin cité plus haut.

\vs
\noindent{\bf N{\oe}uds moyens doubles.} Continuons maintenant en examinant ce qui se passe dans le cas où les deux couples de n{\oe}uds moyens sont doubles, c'est-à-dire quand la souche $A$ est dégagée des couples de branches. Nous sommes alors dans la situation suivante~:~$C=G$, $D=F$ et $C,D$ sont situés de part et d'autre du point $A$ qui est alors le milieu de la pièce $CD$. L'étude proprement dite de cette configuration commence à la page 7, ligne 26 et va se poursuivre jusqu'à la page 10, où Desargues synthétisera les résultats obtenus en généralisant complètement la notion d'involution, avant de passer au \emph{théorème de la ramée,} quittant ainsi définitivement la géométrie en une dimension.

Il commence par rappeler que cette situation de trois couples de n{\oe}uds ne donne en fait que deux couples de points en involution.  Comme l'on a  affaire à un arbre avec deux couples de n{\oe}uds moyens et une souche dégagée, les identités suivantes sont vérifiées~:~
\[
AB.AH=AC.AG=AD.AF=AC^2=AG^2=AD^2=AF^2.
\]
Comme cela a déjà été démontré, nous pouvons en tirer\footnote{p. 7, l. 34}
\[
\frac{AH}{AG}=\frac{AG}{AB}=\frac{HG}{BG}.
\]
Cependant ici $HG=HC$ et $BG=BC$ et donc le dernier rapport est aussi égal à\footnote{p. 7, l. 36}
\[
\sqrt{\frac{HG.HC}{BG.BC}}.
\]
Étant dans une situation d'involution, on peut dire que
\[
\frac{HG.HC}{BG.BC}=\frac{HF.HD}{BF.BD},
\]
et comme $C=G, D=F$, on a aussi $BF=BD, HF=HD$ et $BG=BC$, il s'ensuit\footnote{p. 7, l. 43} que
\[
\frac{HG}{BG}=\frac{HF}{BF}
\]
ce qui, combiné avec l'identité démontrée plus haut, donne finalement\footnote{p. 7, l. 47}
\[
\frac{AH}{AG}=\frac{AG}{AB}=\frac{HG}{BG}=\frac{HF}{BF}.
\]
Desargues donne alors une interprétation géométrique de la situation en termes de trois segments consécutifs~:~$FB, BG, GH$ sont trois pièces consécutives telles que l'un des bouts, $GH$, est à sa mitoyenne $BG$, comme la somme des trois $HF$, est à l'autre bout $BF$\footnote{p. 7, l. 50--55}. On reconnait là une configuration de \emph{division harmonique} et il peut apparaître curieux que Desargues ne le mentionne pas. Il note cependant que cette configuration est affaire d'une paire de couples de points, ici $(B,H)$ et $(F,G)$. Mais comme le mot {\it couple} est déjà utilisé, il va employer une autre terminologie~:~dans le cas d'une involution de 4 points seulement avec deux couples de n{\oe}uds moyens $CG, DF$ et un couple de n{\oe}uds extrêmes $BH$, il déclare les deux n{\oe}uds moyens (mais non couplés) $G,F$ \emph{correspondants entre eux\footnote{p.8, l. 10},} de même que les deux n{\oe}uds extrêmes (qui eux sont couplés) sont déclarés correspondants entre eux. Cela n'est jamais que la mise au jour du fait qu'une division harmonique $F, B, G, H$ ne dépend en fait que de la paire de paires $\{B,H\}, \{F,G\}$, ce qui clarifie complètement la combinatoire de cette configuration qui est bien vue par Desargues comme une \emph{relation} entre points. 

De la ligne 10 à la ligne 39  de la page 8, il répète de plusieurs manières différentes ce qu'il comprend par les mots {\it involution de quatre poincts seulement.} Il entend par là soit la situation décrite jusqu'ici, à savoir une involution avec deux couples de n{\oe}uds moyens doubles et un couple de n{\oe}uds extrêmes, donnant en fait quatre points divisant harmoniquement la droite, soit encore le cas où le n{\oe}ud extrême intérieur se confond avec la souche tandis que l'extérieur part à l'infini, donnant ainsi la configuration que l'on pourrait noter $F,B,G,\infty$, avec $B$ milieu de $FG$. 

À partir de la ligne 40, il revient sur l'étude de l'involution à n{\oe}uds moyens doubles et va montrer comment, en changeant l'involution (et la souche de l'arbre correspondant) on va pouvoir transformer des n{\oe}uds moyens doubles en n{\oe}uds extrêmes et {\it vice versa.} Cela renforce encore ce qui a déjà été démontré sur le fait qu'une division harmonique $F, B, G, H$ ne dépend en fait que de la paire de paires $\{B,H\}, \{F,G\}$, mais surtout illustre notre thèse sur le cheminement de Desargues concernant la genèse de la notion d'involution. 

Prenons en effet le milieu $L$ du segment $HB$. Comme il ressort de ce qui précède,
\[
\frac{BF}{BG}=\frac{HF}{HG}
\]
et donc les rapports des carrés sont égaux eux aussi~:~
\[
\frac{FB.FB}{GB.GB}=\frac{FH.FH}{GH.GH},
\]
ce qu'il interprète\footnote{p. 8, l. 50} comme le fait que $FG, BB, HH$ forment une involution dont $F$ et $G$ sont n{\oe}uds extrêmes et $B,H$ n{\oe}uds moyens doubles. Il remarque ensuite que les mêmes démonstrations faites jusqu'ici impliquent que
\[
\frac{FB.FB}{GB.GB}=\frac{FH.FH}{GH.GH}=\frac{LF}{LG},
\]
ce qui permet alors de dire que $L$ est souche de l'arbre $(L, GF, BB, HH)$.


Cela achève en quelque sorte sa description de la double symétrie de la division harmonique $\{B,H\}, \{F,G\}$ évoquée plus haut. Mais surtout, cela montre que l'on peut considérer une division harmonique \emph{de deux manières} comme une involution, suivant qu'on prend pour souche de l'arbre le milieu $A$ de $FG$ ou le milieu $L$ de $BH$, obtenant ce qu'il appelle deux souches \emph{réciproques entre elles,} ou encore que l'on considère $B,H$ comme n{\oe}uds moyens ou extrêmes (et $F,G$ comme n{\oe}uds extrêmes ou moyens).

De la ligne 8 de la page 9 jusqu'au milieu de la page 10, Desargues reprend son étude des n{\oe}uds moyens doubles. Il commence par affirmer que 
\[
\frac{BF}{BA}=\frac{BH}{BG},
\]
ce qui se déduit en effet de l'égalité $AG.AG=AB.AH$ par les méthodes qu'il a déjà maintes fois utilisées. Transformant cette égalité de rapports en égalité de rectangles on obtient $BF.BG=BA.BH$, ce qui lui fait dire que $B$ peut être considéré comme souche\footnote{«~Pour souche~» : p. 9, l. 11} de l'arbre $(B; FG, AH)$. Ici encore, d'une même division harmonique et du choix d'une souche pour la voir comme une involution, on peut déduire plusieurs involutions. Il va exploiter à fond cette possibilité combinatoire et réutiliser les méthodes de preuve qu'il a mises au point précédemment. Nous donnons la suite des identités qu'il obtient sans plus de précision. Du fait que $B$ est souche de l'arbre $(B; FG, AH)$ on tire
\[
\frac{BF}{BG}=\frac{FA.FH}{GA.GH}
\]
terme égal à $FH/GH$ puisque $FA=GA$, et qu'on peut finalement écrire comme la racine carrée d'un rapport de carrés\footnote{p. 9, l. 15}~:~
\[
\frac{BF}{BG}=\sqrt{\frac{FH.FH}{GH.GH}}.
\]
Comme $B$ et $H$ jouent des rôles symétriques, on peut déduire de l'identité
\[
\frac{BF}{BA}=\frac{BH}{BG}
\]
que tout aussi bien
\[
\frac{HF}{HA}=\frac{HB}{HG}
\]
soit encore\footnote{p.9, l. 17}
\[
\frac{HG}{HB}=\frac{HA}{HF},
\]
ce qui donne l'identité de rectangles $HG.HF=HB.HA$ disant que $H$ est \emph{pour souche} de l'arbre $(H; GF, BA)$. Répétant le même argument que précédemment avec $H$ comme souche on obtient derechef
\[
\frac{HF}{HG}=\sqrt{\frac{BF.BF}{BG.BG}},
\]
puis que\footnote{p. 9, l. 26}
\[
\frac{BF}{BG}=\frac{FA.FB}{GA.GB}\;\mbox{et}\;\frac{HF}{HG}=\frac{FA.FH}{GA.GH}.
\]
Arborisant de plus belle, il en tire\footnote{p.9, l. 29}
\[
\frac{FA.FB}{GA.GB}=\frac{FA.FH}{GA.GH}.
\]
Comme $B$ est souche de $(B; AH, FG)$ on a également
\[
\frac{BH}{BA}=\frac{HF.HG}{AF.AG}
\]
et comme $AF=AG$ ce dernier rapport est aussi égal à $HF.HG/(AG.AG)$. Le même raisonnement avec $H$ pour souche donne
\[
\frac{HA}{HB}=\frac{AF.AG}{BF.BG}=\frac{AG.AG}{BF.BG}.
\]
Réutilisant sans le dire son argument usuel, à partir cette fois-ci de l'identité de rectangles $AF.AG=AF^2=AB.AH$, il tire\footnote{p.9, l. 35}
\[
\frac{FH}{FA}=\frac{BF}{BA}=\frac{BH}{BG},
\]
d'où suit d'abondant que\footnote{p. 9, l. 37}
\[
FH.BA=FA.FB\;\mbox{et}\; HF.BG=FA.BH.
\]
De l'identité
\[
\frac{FH}{FA}=\frac{BF}{BA}
\]
on tire, par différences, 
\[
\frac{AB}{AG}=\frac{AF}{AH}=\frac{BF}{FH}
\]
et donc\footnote{p.9, l. 39}
\[
\frac{GA}{BF}=\frac{HA}{HF}
\]
ainsi que les identités de rectangles qui s'en déduisent. De même, en usant de la souche $H$, obtient-on
\[
\frac{HA}{HF}=\frac{HG}{HB}.
\]
Comme $A$ est milieu de $GF$, il s'ensuit\footnote{p. 9, l. 42}
\[
\frac{BF}{BH}=\frac{2FA}{2GH}=\frac{FG}{2GH}
\]
et ainsi
\[
\frac{FB}{FG}=\frac{BH}{2HG}
\]
et donc\footnote{p. 9, l. 44}
\[
\frac{FB}{FA}=\frac{2FB}{FG}=\frac{2BH}{2GH}=\frac{BH}{GH}
\]
et les identités de rectangles
\[
FB.HG=FA.HB\;\mbox{et}\; FB.HA=FA.HF
\]
qui s'ensuivent. Desargues, tel un géant jonglant avec des arbres, énonce alors toute une série d'identités que nous reproduisons ici pour le \emph{divertissement\footnote{dernière ligne de la page 9}} du lecteur~:~
\[
\frac{HB.HB}{HB.HA}=\frac{BA.BH}{AB.AH}=\frac{BG.BF}{AG.AG}=\frac{BG.BF}{AF.AF}=\frac{HB}{HA},
\]
puis
\[
\frac{BG}{BA}\frac{BF}{AH}=\frac{HB}{HA}
\]
et comme
\[
\frac{HB}{HA}=\frac{BG.BF}{AG.AF}=\frac{GB}{GA}\frac{FB}{FA}
\]
cela donne\footnote{p. 9, l. 59}
\[
\frac{BG}{BA}\frac{BF}{AH}=\frac{GB}{GA}\frac{FB}{FA}=\frac{HB}{HA}.
\]
Après une phrase mise comme une respiration au bas de la page 9, l'auteur continue, comme enivré de sa virtuosité. Le premier alinéa de la page 10 est ainsi tourné~:~«~Davantage puis que HB, est à HG, comme FB, est à FA, \& que la raison est double qui est composée des raisons de FG, à FB, \& de FB, à FA, c'est à dire la raison de FG, à FA.~» Que signifie ici le mot «~double~»? Desargues dit d'abord que
\[
\frac{HB}{HG}=\frac{FB}{FA},
\]
puis que (c'est la raison composée)
\[
\frac{FG}{FB}\frac{FB}{FA}=\frac{FG}{FA},
\]
et ce dernier rapport vaut 2, puisque $A$ est milieu de $FG$. C'est la raison double dont il est question dans le texte. De là il tire\footnote{p. 10, l. 4}
\[
\frac{FG}{FB}\frac{FB}{FA}=\frac{FG}{FB}\frac{HB}{HG}=2=\frac{FG}{HG}\frac{HB}{FB}.
\]
Avec le même argument, partant de 
\[
\frac{BH}{BG}=\frac{FH}{FA},
\]
il obtient\footnote{p. 10, l. 10}~:~
\[
\frac{FG}{BG}\frac{BH}{FH}=2=\frac{HB}{HF}\frac{GF}{GB}.
\]
Ici s'achève cette chevauchée dans ce qu'il convient maintenant d'appeler une forêt. 

Il va alors tirer quelques conséquences des propriétés qu'il vient de démontrer pour savoir si quatre points sont en involution\footnote{p.10, l. 14}. 

Le premier exemple\footnote{p. 10, l. 16} s'énonce ainsi~:~si $AB/AC=AC/AH$ ($AB, AC$ et $AH$ sont dites «~entre elles continuellement proportionnelles~») et si $AF=AC$, alors $H,C,B,F$ sont en involution. Ensuite\footnote{p. 10, l. 19}, si $BH/BG=BF/BA$ (les quatres segments $BH,BG, BF, BA$ sont deux à deux proportionnels) et si $A$ est milieu de $FG$, alors $H,G,B,F$ sont en involution. Enfin\footnote{p. 10, l. 23}, si $HG/HB=HA/HF$ et si $A$ est milieu de $FG$ alors $H, G, B, F$ sont en involution. 

Ici s'achève ce long et difficile passage du {\it Brouillon-Project,} à la fin duquel Desargues fait la synthèse citée plus haut. Il passe alors à la géométrie plane proprement dite, en démontrant le {\it théorème de la ramée,} qui n'est jamais que l'énonciation de l'invariance de l'involution par projection centrale. 

\section{Arbres et involutions~:~une version utilisant le langage algébro-géométrique}
Considérons un corps commutatif $\KK$ (de caractéristique différente de $2$) et $\LL$ une droite projective sur le corps $\KK$. Munissons $\LL$ d'un repère projectif $\{0,1,\infty\}$ permettant d'identifier $\LL$ à $\KK\cup\{\infty\}$. Soit $a\in\LL$ un point, dit {\it souche,} et $b,h;c,g;d,f$ six points distincts deux à deux et différents de $a$, accouplés deux à deux. 
\begin{definition} On dit que la configuration $a;b,h;c,g;d,f$ forme un \emph{arbre algébrique} si les identités suivantes sont satisfaites~:~
\[
(b-a)(h-a)=(c-a)(g-a)=(d-a)(f-a).
\]
\end{definition}
Suivons le cheminenent de Desargues et transformons ces identités de rectangles en identités de rapports de rectangles ne faisant plus apparaître la souche $a$. Prenons deux couples quelconques, par exemple $c,g$ et $d,f$. De l'identité des rectangles on tire deux identités de rapports~:~
\[
\frac{c-a}{d-a}=\frac{f-a}{g-a}\;\mbox{et}\;\frac{c-a}{f-a}=\frac{d-a}{g-a}.
\]
Nommons $\alpha$ le premier rapport et $\alpha'$ le second, de sorte que
\[
\alpha=\frac{c-a}{d-a}=\frac{f-a}{g-a}=\frac{(c-a)-(f-a)}{(d-a)-(g-a)}=\frac{c-f}{d-g},
\]
et que
\[
\alpha'=\frac{c-a}{f-a}=\frac{d-a}{g-a}=\frac{(c-a)-(d-a)}{(f-a)-(g-a)}=\frac{c-d}{f-g}.
\]
Les opérations effectuées ici sont licites car les six points sont disctincts deux à deux. En multipliant $\alpha$ et $\alpha'$ on obtient
\[
\alpha\alpha'=\frac{c-a}{d-a}\frac{d-a}{g-a}=\frac{c-a}{g-a}=\frac{(c-f)(c-d)}{(d-g)(f-g)},
\]
soit encore
\begin{equation}\label{acag}
\frac{c-a}{g-a}=\frac{(d-c)(f-c)}{(f-g)(d-g)},
\end{equation}
qui n'est autre que l'identité énoncée par Desargues à la ligne 15, page 4, du {\it Brouillon Project.} En utilisant l'autre identité de rectangles, à savoir $(b-a)(h-a)=(c-a)(g-a)$, on obtient de même
\[
\beta=\frac{c-a}{b-a}=\frac{h-a}{g-a}\;\mbox{et}\; \beta'=\frac{c-a}{h-a}=\frac{b-a}{g-a}
\]
d'où l'on tire
\[
\beta=\frac{c-h}{b-g}\;\mbox{et}\; \beta'=\frac{c-b}{h-g}.
\]
En multipliant $\beta$ par $\beta'$ on obtient comme ci-dessus
\[
\frac{c-a}{g-a}=\frac{(c-h)(c-b)}{(b-g)(h-g)}
\]
ce qui, combiné avec l'identité \ref{acag} entraîne l'égalité de rapports~:~
\[
\frac{(h-c)(b-c)}{(b-g)(h-g)}=\frac{(d-c)(f-c)}{(f-g)(d-g)}.
\]
En faisant tourner la combinatoire gémellaire, on obtient deux autres égalités de rapports de rectangles, nous permettant de donner la définition suivante~:~
\begin{definition} La configuration $b,h;c,g;d,f$ forme une \emph{involution algébrique} si les six points $b,h,c,g,d,f$ sont distincts deux à deux  et si les trois identités suivantes sont satisfaites~:~
\begin{equation}\label{IA1}
\frac{(d-g)(f-g)}{(d-c)(f-c)}=\frac{(b-g)(h-g)}{(b-c)(h-c)},
\end{equation}
\begin{equation}\label{IA2}
\frac{(c-f)(g-f)}{(c-d)(g-d)}=\frac{(b-f)(h-f)}{(b-d)(h-d)}
\end{equation}
et
\begin{equation}\label{IA3}
\frac{(c-h)(g-h)}{(c-b)(g-b)}=\frac{(d-h)(f-h)}{(d-b)(g-b)}.
\end{equation}
\end{definition}
La manière dont Desargues organise ses rapports rend le lien avec la géométrie projective moderne un peu indirect. Prenons l'identité \ref{IA1} et réorganisons-en les termes, elle devient
\[
\frac{(b-g)(d-c)}{(b-c)(d-g)}=\frac{(f-g)(h-c)}{(f-c)(h-g)},
\]
ce que l'on peut écrire encore
\[
\frac{\left(\displaystyle{\frac{c-d}{c-b}}\right)}{\left(\displaystyle{\frac{g-d}{g-b}}\right)}=
\frac{\left(\displaystyle{\frac{c-h}{c-f}}\right)}{\left(\displaystyle{\frac{g-h}{g-f}}\right)},
\]
ce qui n'est autre que l'identité de birapports
\[
[c,g;d,b]=[c,g;h,f].
\]
Sous cette forme condensée, les autres identités se déduisent facilement, et l'on peut donc donner une nouvelle version de l'involution algébrique~:~
\begin{definition} La configuration $b,h;c,g;d,f$ forme une \emph{involution algébrique} si les six points $b,h,c,g,d,f$ sont distincts deux à deux  et si les trois identités suivantes sont satisfaites~:~
\[
[b,h;d,c]=[b,h;g,f], [c,g;d,b]=[c,g;h,f]\;\mbox{et}\; [d,f;c,b]=[d,f;h,g].
\]
\end{definition} 
De cette nouvelle définition il ressort que la notion d'involution algébrique est indépendante du choix du repère projectif sur la droite $\LL$, puisque le birapport est invariant par homographie. Du théorème suivant découlera alors que la notion d'arbre algébrique elle-même est indépendante du choix du repère.
\begin{theoreme}\label{arbre-invol-alg} Trois couples de six points distincts $b,h;c,g;d,f$ forment une involution algébrique si et seulement s'il existe un septième point $a$ faisant office de souche à l'arbre $a;b,h;c,g;d,f$.
\end{theoreme}
\begin{demonstration} Le sens indirect a été démontré plus haut et, pour démontrer le sens direct, nous allons supposer que $(b,h;c,g;d,f)$ forme une involution algébrique. Quitte à changer de repère projectif en particulier en changeant le point à l'infini, on peut supposer que pour deux couples $c,g;d,f$ au moins, on a $d-g\neq c-f$. On procède alors comme Desargues et l'on note $a$ l'unique solution de l'équation linéaire
\[
\frac{g-a}{f-a}=\frac{d-g}{c-f}.
\]
Nous en tirons
\begin{equation}\label{gafa}
\frac{g-a}{f-a}=\frac{d-g}{c-f}=\frac{(g-a)+(d-g)}{(f-a)+(c-f)}=\frac{d-a}{c-a},
\end{equation}
ce qui peut encore s'écrire $(g-a)(c-a)=(d-a)(f-a)$ et il est facile de voir que $a$ est nécessairement différent de $g,f,d,c$. Continuons comme Desargues et montrons par l'absurde que $(b-a)(h-a)=(g-a)(c-a)$. Supposons que cela ne soit pas le cas et considérons $h'$ l'unique solution de l'équation linéaire $(b-a)(x-a)=(g-a)(c-a)$. La configuration $(a;b,h';c,g;d,f)$ forme donc un arbre algébrique et par construction $h\neq h'$. En reprenant, à partir de l'identité $(b-a)(h'-a)=(g-a)(c-a)$ le raisonnement démontrant qu'un arbre algébrique est une involution, on montre que
\[
\frac{c-a}{b-a}=\frac{h'-a}{g-a}=\frac{c-h'}{b-g}
\]
et que
\[
\frac{c-a}{h'-a}=\frac{b-a}{g-a}=\frac{c-b}{h'-g}.
\]
Puisque
\[
\frac{c-a}{g-a}=\frac{c-a}{b-a}\frac{b-a}{g-a}=\frac{c-h'}{b-g}\frac{c-b}{h'-g},
\]
on a
\begin{equation}\label{caga}
\frac{c-a}{g-a}=\frac{c-b}{b-g}\frac{c-h'}{h'-g}\neq \frac{c-b}{b-g}\frac{c-h}{h-g}.
\end{equation}
Mais de l'identité \ref{gafa} on tire
\[
\frac{g-a}{f-a}=\frac{d-a}{c-a},
\]
ce qui permet de démontrer que
\[
\frac{c-a}{g-a}=\frac{(c-f)(c-d)}{(d-g)(f-g)}.
\]
Comme $(b,h;c,g;d,f)$ est une involution algébrique, l'identité \ref{IA1} entraine alors que
\[
\frac{c-a}{g-a}=\frac{(c-b)(c-h)}{(b-g)(h-g)},
\]
en contradiction avec l'inégalité \ref{caga}. Le théorème \ref{arbre-invol-alg} est démontré.
\end{demonstration}
Nous voici maintenant un peu plus à l'aise pour explorer la notion d'involution. Soit $(b,h;c,g;d,f)$ une involution, et soit $a$ la souche en faisant un arbre. Quitte à changer de repère projectif, nous supposerons que $a=0$. Les identités arboricoles deviennent alors
\[
bh=cg=df.
\]
Notons $\alpha$ la valeur commune à ces trois produits et notons $\phi$ l'homographie de $\LL$ dans elle-même définie par
\[
\phi(z)=\frac{\alpha}{z}\;\mbox{ou encore}\; \phi([z:t])=[\alpha t:z].
\]
La transformation $\phi$ est une \emph{involution} au sens moderne du terme et les identités ci-dessus disent simplement que $\phi$ échange $b$ et $h$, $c$ et $g$ ainsi que $d$ et $f$. Elle envoie de plus la souche $a$ sur l'infini.
\begin{proposition} Une configuration $(b,h;c,g;d,f)$ est une involution algébrique si et seulement s'il existe une homographie involutive $\phi$ telle que $\phi(b)=h, \phi(c)=g$ et $\phi(d)=f$.
\end{proposition}
Une homographie étant entièrement déterminée par ses valeurs sur trois points distincts, la donnée de l'involution algébrique et celle de l'homographie involutive associée sont équivalentes. On en déduit immédiatement le résultat suivant, cas particulier de l'énoncé du {\it Brouillon Project} cité plus haut~:~
\begin{proposition}\label{involution-generale} Soit $(b,h;c,g;d,f)$ une involution algébrique et supposons-nous donné un autre couple de points $k,l$ tel que $(b,h;c,g;k,l)$ forme une involution. Alors $(c,g;d,f;k,l)$ forme également une involution.
\end{proposition}
Jusqu'ici nous avons supposé les points des divers couples distincts deux à deux. Une homographie involutive est toujours conjuguée à une transformation telle que $\phi$ ci-dessus, avec $\alpha\in\KK^*$. Elle admet deux points fixes si $\alpha$ est un carré dans $\KK$, et aucun sinon. Rappelons que dans le premier cas on dit que $\phi$ est \emph{hyperbolique} et que dans le second elle est \emph{elliptique,} terminologie qui prendra son sens par la suite. Il apparaît dès lors une première dichotomie dans la notion d'involution arguésienne, suivant la nature de l'homographie involutive qu'elle détermine. 
\begin{definition} Une involution (au sens de Desargues, ou bien au sens algébrique) est dite elliptique si l'homographie involutive qu'elle détermine est sans point fixe. Elle est dite hyperbolique si cette homographie a deux points fixes.
\end{definition}
\subsection{Involutions hyperboliques}
Supposons que $(b,h;c,g;d,f)$ est une involution hyperbolique et notons $\phi(z)=\alpha/z$ son homographie involutive associée. L'équation $x^2=\alpha$ admet deux solutions, donnant les deux points fixes $k$ et $l$ de $\phi$ (en coordonnées, $k,l=\pm \sqrt{\alpha}$). Rappelons la proposition suivante~:~
\begin{proposition} Étant donnés deux points distincts $k,l$ de la droite $\LL$, il existe une unique homographie involutive $\phi$ admettant $k$ et $l$ pour points fixes. De plus, pour tout point $b$ d'image $h$, on a
\[
[k,l;b,h]=-1.
\]
\end{proposition}
Nous pouvons donc étendre la notion d'involution arguésienne au cas où deux des couples $c,g$ et $d,f$ sont des couples de points confondus. Dans ce cas les rectangles $AC.AG$ et $AD.AF$ correspondants sont en fait des carrés~:~$c,g$ et $d,f$ sont alors ce que Desargues appelle des \emph{n{\oe}uds moyens.} Mieux, comme ici l'involution est hyperbolique, on a nécessairement $c=g=k$ et $d=f=l$, c'est-à-dire que dans le cas d'une involution hyperbolique, l'involution possède deux couples de \emph{n{\oe}uds moyens doubles\footnote{Et réciproquement.}}. Si $b$ est un tierce point tel que $h=\phi(b)\neq b$, la configuration $(b,h;k,k;l,l)$ forme donc une involution de quatre points seulement, et la proposition ci-dessus nous informe que
\[
[k,l;b,h]=-1,
\]
c'est-à-dire que l'involution de quatre points seulement $(b,h;k,l)$ est une \emph{division harmonique} de la droite $\LL$. Dans le cas où le corps de base $\KK$ est le corps des nombres réels, dire que $\phi$ est hyperbolique revient à dire que $\alpha$ est positif, et l'on retrouve alors le fait que dans ce cas la souche est dégagée d'entre les couples~:~$c=g=k$ est d'un côté de la souche $a$ (en coordonnées, $a=0$), $d=f=l$ de l'autre, et les deux points $b,h$ sont du même côté de la souche, l'un $b$ est plus proche de $a$ que ne l'est $c$ par exemple, l'autre $h$ en est plus éloigné. Remarquons que cela permet encore d'étendre la notion d'involution, en décrétant que $(a,\infty;k,k;l,l)$ forme une involution si $k^2=l^2=\alpha$ et si $a$ est le milieu de du segment $[l,k]$. 

Faisons encore une remarque~:~soit $(f,g;b,h)$ une division harmonique de la droite. Alors il existe une unique involution $\sigma$ telle que $b,h$ soient ses points fixes et $\sigma(f)=g$. De même il existe une unique involution $\sigma'$ telle que $f,g$ soient ses points fixes et $\sigma'(b)=h$. Fixons une structure affine en fixant un point $\infty$ comme point à l'infini sur la droite. Si l'on pose $a=\sigma(\infty)$, on dit que $a$ est le point central de l'involution $\sigma$ et elle est déterminée par l'équation $(m-a)(\sigma(m)-a)=\alpha$. Si l'on pose $\ell=\sigma'(\infty)$, alors $\ell$ est point central de l'autre involution $\sigma'$ déterminée par la division harmonique $(f,g;b,h)$, et l'on retrouve ce que Desargues nomme les \emph{souches réciproques entre elles} de l'involution de quatre points seulement $(f,g;b,h)$. Suivant que l'on considère que $f,g$ sont n{\oe}uds moyens ou extrêmes, on obtient d'une telle situation en division harmonique deux manières d'involution  arguésienne. Dit autrement, la notion d'involution arguésienne, même de quatre points seulement, est plus riche que celle de division harmonique.
\subsection{Involutions elliptiques}
Supposons que $(b,h;c,g;d,f)$ est une involution elliptique et notons $\phi(z)=\alpha/z$ son homographie involutive associée. L'équation $x^2=\alpha$ est sans solution. Il n'existe donc pas de point fixe ni \textit{a fortiori,} de n{\oe}uds moyens dans ce cas. Rappelons cependant que Desargues définit la notion de n{\oe}ud moyen en termes de longueurs, c'est-à-dire d'égalités de quantités positives ou, plus généralement, de quantités qui sont des carrés. Pour simplifier les choses, 
plaçons-nous dans le cas où $\KK=\RR$. Le terme $\alpha$ est dans notre cas strictement négatif. Il existe $g>0$ et $c<0$ tels que $|g|=|c|$ et $gc=\alpha$. Dit autrement le couple $c,g$ (en coordonnées, le couple $\sqrt{-\alpha},-\sqrt{-\alpha}$) est un couple de n{\oe}uds de l'involution déterminée par $\phi$, et comme $|g|=|c|$, Desargues les considère comme des n{\oe}uds moyens. Comme $g\neq c$ puisque $c=-g$, il le nomme un couple de \emph{n{\oe}uds moyens simples.} Ils sont situés de part et d'autre de la souche $a$ et celle-ci est donc engagée. Il existe alors un unique autre couple de n{\oe}uds moyens simples, à savoir $d=-\sqrt{-\alpha},f=\sqrt{-\alpha}$ de sorte que l'involution possède deux couples de n{\oe}uds moyens simples, le couple $c,g$ et le couple $d,f$, ne donnant en fait que deux points puisque $c=f$ et $d=g$. 

Pour se tirer de l'embarras ainsi crée par la coïncidence de deux points pris dans deux couples différents, et inclure ce cas dans le cadre des involutions de quatre points seulement, comme dans le cas hyperbolique, Desargues écrit la chose suivante, lignes 31 et suivantes de la page 10 du {\it Brouillon Project~:~}«~Quand en un tronc droict, trois couples de n{\oe}uds extremes DF, CG, BH, sont en involution entre eux, que deux autres couples de n{\oe}uds, moyens, unis, doubles ou simples, PQ,XY, sont une involution de quatre points avec chacune des deux quelconques couples CG,\& BH, de ces trois couples de n{\oe}uds extremes. Ces deux mesmes n{\oe}uds moyens PQ, XY,
sont encore une involution aussi de quatre points avec la troisiéme de ces couples de n{\oe}uds extremes DF.~»

Le problème est que la notion de n{\oe}ud moyen simple n'est pas de nature projective. Plus précisément~:~que $f$ soit elliptique ou hyperbolique, cela est un fait projectif. Si elle est hyperbolique, elle admet deux points fixes et cela encore est un fait projectif~:~ces deux points fixes, donnant les n{\oe}uds moyens doubles, sont des invariants de conjugaison homographique. Mais si l'involution est elliptique, déterminer des n{\oe}uds moyens simples suppose que l'on a fixé le point à l'infini, c'est-à-dire une structure affine sur la droite projective privée de ce point. Voilà sans doute une des raisons pour lesquelles Desargues bute sur ce point, raison que bien évidemment il ne peut entendre. 

\subsection{Involutions «~paraboliques~»}
Rappelons qu'une homographie $\phi$ est dite {\it parabolique} si elle possède un unique point fixe. Il est bien connu que si $\phi$ est involutive, elle ne peut être parabolique\footnote{Sauf sur un corps de caractéristique 2, où la translation par 1 donne une involution, mais on ne peut guère soupçonner Desargues d'avoir eu en tête une telle horreur.}. On trouve cependant dans l'article de Lenger \cite{lenger} l'idée que Desargues aurait eu l'intuition du «~troisième cas possible d'involution~», idée reprise dans \cite{taton}, p. 157. Cette erreur est dûe à une interprétation erronée d'un passage qui commence ligne 30 de la page 21 du {\it Brouillon Project~:~}«~Cependant on remarquera qu'entre les deux espèces de conformation d'arbre il y en a une troisiesme en laquelle de chaque couple de n{\oe}uds, tousiours un est uny à la souche, ou l'entendement demeure court de mesme qu'en plusieurs autres circonstances, \& cette espece de conformation d'arbre est mytoyenne entre les autres à souche engagée \& souche dégagée.~» On remarque tout d'abord qu'ici Desargues revient à la notion d'arbre et ne parle aucunement d'involution. Revenons à la notation algébrique moderne en prenant le point $a=0$ pour souche et considérons une involution $\phi$ de point central $0$. On sait alors qu'il existe un nombre réel non-nul $k$ tel que $z\phi(z)=k$. Si $k<0$ l'involution est elliptique et sans point fixe, la souche est engagée. Si $k>0$, elle est hyperbolique et a deux points fixes, la souche est dégagée. Le cas dont parle Desargues, ({\it mytoyen entre les autres}) est donc celui où $k=0$, ce qui définit une transformation qui n'est plus du tout une homographie mais celle qui envoie tout point $z$ (sauf l'infini) sur la souche $a$ ({\it de chaque couple de n{\oe}uds, tousiours un est uny à la souche}). C'est une transformation singulière et Desargues se garde bien d'inclure ce cas au rang d'une situation d'involution même s'il apparaît géométriquement comme cas dégénéré~:~un peu plus bas, Desargues considère ce qu'il advient quand une {\it traversale} est tangente à la conique. Précisons brièvement de quoi il s'agit. Si $\cC$ est une conique et si $F$ est un point du plan, sa {\it traversale} n'est jamais que sa polaire eu égard à $\cC$. Si $F$ est sur la conique, sa polaire est la tangente à la conique en ce point. Sur le plan vectoriel définissant cette tangente, la (classe de la) forme quadratique définissant la conique dégénère et ne peut donc plus définir une involution, mais une transformation singulière ou arbre {\it mytoyen} ainsi que défini ci-dessus. Desargues écrit à ce sujet~:~«~cet arbre est de l'espèce mitoyenne, dont l'entendement ne peut comprendre comment sont les proprietez que le raisonnement luy en fait conclure.~» Il ne découvre donc pas un hypothétique cas d'involution parabolique, mais bute sur une véritable difficulté conceptuelle. Notons que le terme d'involution parabolique apparaît aussi chez Zacharias dans \cite{zacharias}, mais que cet auteur prend bien garde de mettre ces mots entre guillemets. 
\begin{center}
{\Large\bf Partie II 

\vs
M. de Beaugrand, lecteur des dix premières pages du {\it Brouillon}}
\end{center}

Nous voudrions dans cette partie montrer comment le {\it Brouillon  Project} de Desargues a été reçu par certains de ses contemporains, en analysant la lecture qu'en fait Jean de Beaugrand. La source que nous utilisons est le petit fascicule imprimé conservé à la Bibliothèque nationale de France, au département Littérature et art, sous la cote V-12591. Il est intitulé {\it Advis Charitables sur les diverses {\oe}uvres, et feuilles volantes du S${}^r$ Girard Desargues, Lyonois.} Ce fascicule comporte trois parties dont la première porte le titre de {\it Lettre de Monsieur de Beaugrand Secretaire du Roy, sur le suject des feuilles intitulees. Brouillon Project d'une atteinte aux evenements des rencontres du Cone avec un plan, \& aux evenements des contrarietez, d'entre les actions des puissances.} Beaugrand y examine donc le {\it Brouillon Project,} y compris son appendice portant sur le problème de la Géostatique, sujet de controverse à l'époque ({\it voir} par exemple \textit{La vie de M.~Descartes} d'Adrien Baillet). 
\section{Beaugrand et les {\it Advis Charitables}}
Jean de Beaugrand est né vers 1580 et mort en 1640. Membre de l'académie de Mersenne et secrétaire du Roy en charge de la délivrance des {\it privilèges} autorisant la publication, il a joué un rôle important dans la vie scientifique de la première moitié du dix-septième siècle. Ses controverses avec Descartes ou Desargues, son ton volontiers polémique et sa propension supposée à l'appropriation des idées d'autrui ont largement entâché sa postérité\footnote{Nous renvoyons pour plus de détail à la page {\it wikipedia} qui lui est consacrée, seule source synthétique connue de nous concernant ce personnage.}.

La lettre attaquant Desargues a été publiée à titre posthume, en 1642. C'est la seule des trois parties des {\it Advis Charitables} qui soit signée. Les deux autres portent respectivement sur les travaux de Desargues en perspective et en stéréotomie et sur ses travaux touchant à la gnomonique. Un trait commun à ces textes est qu'ils attaquent frontalement Desargues dans sa prétention à la nouveauté et à l'universalité, en pointant plus particulièrement la terminologie exotique introduite par le lyonnais. On y trouve de nombreuses citations. Celle concluant l'avertissement au lecteur ouvrant les {\it Advis,} tirée du livre de la {\it Genèse,} illustre parfaitement l'esprit de l'ensemble~:~«~\textit{Terra autem erat inanis et vacua, et tenebr{\ae} erant super faciem abyssi\footnote{Cependant la terre était informe et vide, et les ténèbres recouvraient l'abysse.}}~» Dit autrement, les écrits de Desargues sont obscurs et vides de contenu et rien de nouveau ne s'y trouve que l'on ne trouvât chez les auteurs de la tradition, plus particulièrement les antiques pour ce qui concerne la géométrie.
\section{La lecture par Beaugrand des dix premières pages du {\it Brouillon--project}}
Comme pour les autres textes constituant les {\it Advis Charitables,} Beaugrand s'attaque d'abord au style de Desargues et notamment à son vocabulaire novateur. Il use pour cela de traîts d'humour et d'ironie en ayant recours à de nombreuses citations. Il commence par se référer au livre 57  de l'{\it Histoire Romaine} de Dion Cassius, où est décrit comment un certain Marcellus s'opposa au futur empereur Tibère au sujet de l'insertion d'un mot dans le dictionnaire. Il écrit alors\footnote{p.1, l. 13 de la {\it lettre.}} ne pas pouvoir excuser celui qu'il qualifie d'ami\footnote{Beaugrand est resté longtemps en bons termes avec Desargues.} d'avoir non seulement substitué des termes barbares à ceux reçus des anciens, mais d'en avoir en outre introduit de nouveaux, «~entièrement ridicules~». Il énumère alors les nombreux termes introduits par Desargues et conclut \footnote{p.1, l. 22} de manière humoristique en qualifiant ce vocabulaire de «~plus capable de mettre les esprits {\it en involution} ou d'en faire des {\it souches reciproques,} que de leur donner quelque nouvelle lumiere dans les Mathematiques.~»  

Beaugrand l'écrit lui-même, il n'a lu avec soin que les dix premières pages du {\it Brouillon} ce qui, en constituant le tiers, est à ses yeux suffisant pour juger des apports de ce texte. Eût-il lu le reste avec sérieux que son opinion aurait été sans doute changée. Selon lui, ces dix premières pages, celles qui traitent de l'involution, ne sont jamais que la conséquence d'une «~proposition qui est parmy les lemmes du livre de la section determinée dans le septiesme de Pappus~». Il s'agit en fait du premier des lemmes, ou encore de la proposition 22 du livre 7 de la {\it Collection} de Pappus ({\it voir} l'édition de Ver Eecke \cite{pappus}, pp. 530 \& 531, ainsi que celle de Jones \cite{pappus2}, p. 142), que Beaugrand ~«~paraphrase~», comme il le dit lui-même, de la façon suivante~:~




«~Si en la droicte infinie AB, les rectangles ADC, BDE sont esgaux, je dis que BD est à DE, comme le rectangle ABC est au rectangle AE\footnote{Ici, il faut lire AEC},\& qu'aussi AD, est à DC, comme le rectangle BAE, est au rectangle\linebreak
BCE\footnote{p. 2, l. 27.}.~»

 Le lien avec les notions arguésiennes se fait facilement~:~$D$ est souche d'un arbre dont $A,C$ et $B,E$ sont deux couples de n{\oe}uds. Nous retrouvons donc ici ce qui occupe Desargues dans les lignes 7 à 34 de la page 4 du Brouillon--Project. Beaugrand semble ne pas se soucier du problème de la disposition des points, qu'il traitera un peu plus loin, de manière quelque peu cavalière.

\begin{figure}[h]
\centering
\includegraphics[width=7cm]{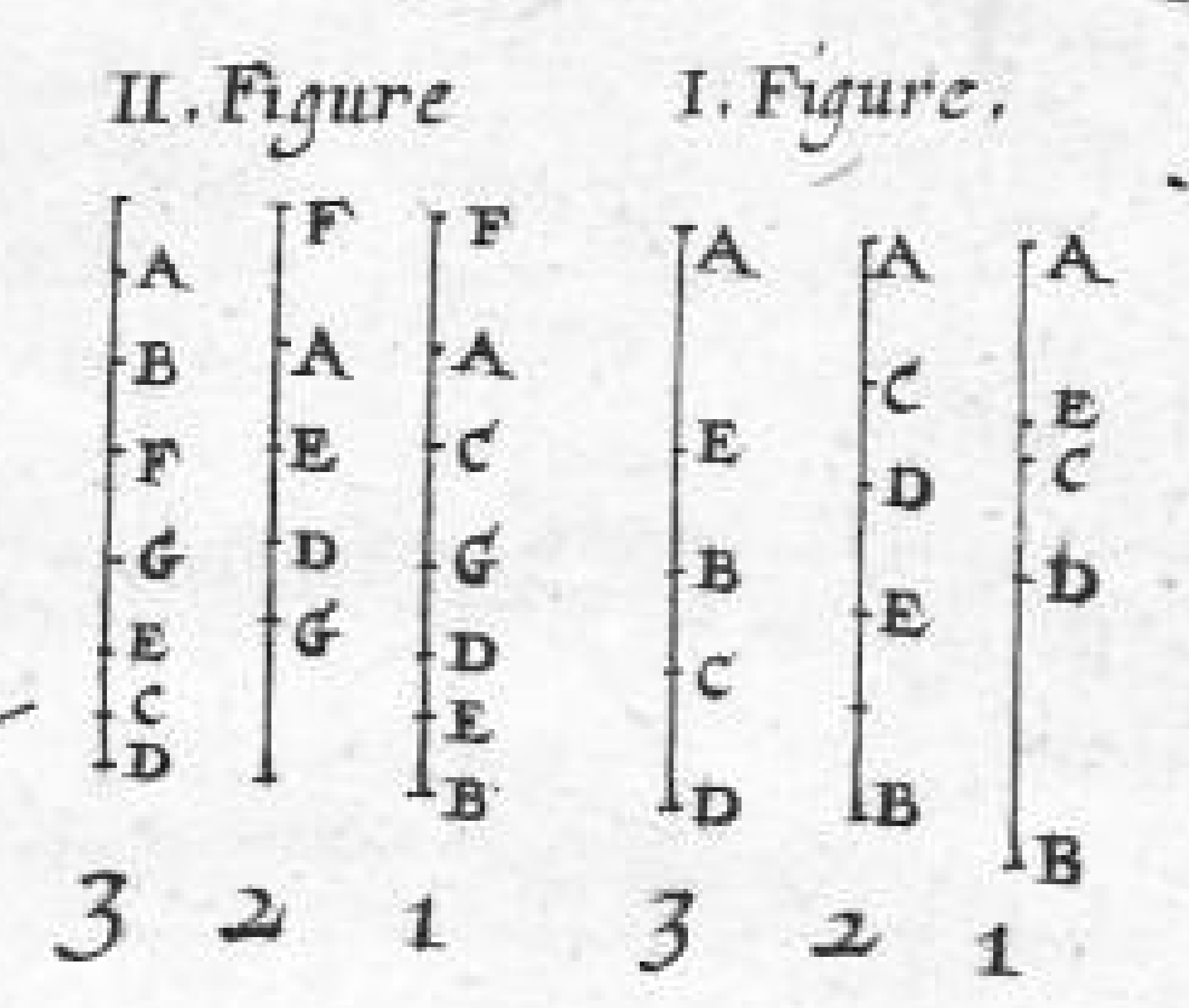}
\caption{Figures I \& II du texte des \textit{Advis charitables} de Beaugrand}
\label{Figures-1-2-Advis}
\end{figure}

La démonstration de Beaugrand suit celle de Pappus. Nous la donnons ici pour la commodité du lecteur. De l'hypothèse $AD.DC=BD.DE$ on tire la proportion $AD/DB=DE/DC$. Par somme ou différence on obtient ainsi
\[
\frac{AD}{DB}=\frac{DE}{DC}=\frac{AD\pm DE}{BD\pm DC}=\frac{AE}{BC}.
\]
De l'autre proportion tirée de l'hypothèse on déduit de même
\[
\frac{AD}{DE}=\frac{BD}{DC}=\frac{AB}{CE}.
\]
On a donc les égalités
\[
\frac{BC}{AE}=\frac{DC}{DE}\;\mbox{et}\;\frac{BD}{DC}=\frac{AB}{CE},
\]
d'où l'on tire
\[
\frac{BD}{DE}=\frac{BD}{CD}\frac{DC}{DE}=\frac{AB}{CE}\frac{BC}{AE}
\]
soit finalement
\[
\frac{BD}{DE}=\frac{AB.BC}{AE.EC},
\]
ce qui est la première égalité cherchée\footnote{Ceci se trouve à la dernière ligne de la page 2 de l'original, avec une faute de typographie, Beaugrand écrivant que $BD/DE=\frac{AD.DC}{AE.EC}$.}. Il réutilise ensuite les deux proportions $AD/DE=BA/CE$ et $DE/DC=AE/BC$ pour déduire derechef que $(AD/DE)(DE/DC)=(BA/CE)(AE/BC)$ soit la deuxième égalité recherchée, à savoir\footnote{Beaugrand utilise une terminologie un peu particulière pour la composition des raisons~:~il écrit ainsi, en bas de la page 2 de l'original~:~«~la raison de BD à DC plus la raison de DC à DE sera esgale à la raison de AB à CE : plus la raison de BC à AE.~» Le mot «~plus~» est ici employé pour désigner la composition des raisons et non, ce qui n'aurait pas de sens de toute façon, leur addition.}
\[
\frac{AD}{DC}=\frac{BA.AE}{BC.CE}.
\]
Après cette démonstration il précise qu'il existe plusieurs dispositions possibles de ces points et nous renvoie à la figure I  ({\it voir} la figure \ref{Figures-1-2-Advis}).

Le premier cas décrit par Beaugrand est celui correpondant à la figure I.2 qui, en termes arguésiens, est celui à souche dégagée (rappelons qu'ici $D$ est la souche). Le second cas est celui à souche engagée, quand $D$ est entre $A$ et $C$ et entre $B$ et $E$. Il ne correspond à aucun des trois cas de la figure I~! Vient ensuite le troisième et, d'après Beaugrand, dernier cas, qui renvoie à la figure I.3 et qui correspond encore une fois à une situation de souche dégagée.

Beaugrand affirme ensuite que la réciproque du lemme énoncé est une évidence, invoquant un raisonnement par contraposée (qu'il nomme «~consécution renversée~», \textit{voir} ligne 17 et suivantes de la page 3 de l'original). Il ne lui reste plus qu'à dire que si l'on introduit un troisième rectangle comme $FD.DG$ égal au deux autres, on obtiendra le même genre de proportion, concluant alors que c'est là «~tout ce qui est contenu dans le tiers de ce grand \textit{broüillon project}\footnote{p. 3, l. 41}~». Il évoque ensuite, dans la première moitié de la page 4, comment comprendre la combinatoire des analogies démontrées et les retrouver presque automatiquement. Beaugrand a ainsi (presque) démontré que $(D; AC,BE,FG)$ est un arbre si et seulement si $(AC, BE, FG)$ est une involution et rejette cette terminologie comme inutile. Il a au passage rapidement réglé le cas des situations particulières, comme pouvant être «~aysém\~et faict par ceux qui ont appris les premiers elem\~es d'Euclide\footnote{p.4, l. 1}~». 

Il poursuit, dans la deuxième moitié de la page 4 et la première de la page 5, par une longue saillie sarcastique sur le vocabulaire arguésien, d'abord sur le mot d'involution, puis sur ceux choisis par Desargues pour remplacer les termes issus du grec que sont l'ellipse, la parabole et l'hyperbole. Il note au passage que la proposition qu'il vient de démontrer de manière courte et synthétique a demandé à Desargues un long développement et «~ample provision \textit{d'arbres, de troncs, de souches, de racines, \&c.}~»

Il passe alors à la démonstration d'une version du théorème d'involution de Desargues, en l'énonçant ainsi~:~«~Si on prend les quatre points K, N, O, V, dans la section Conique, K N G O V F, \& que l'on tire les quatre droictes KN, KO, VN, VO tellement qu'il parte deux droictes \& non plus de chacun de ces poincts, \& puis que l'on tire en telle façon que l'on voudra la droite F ACG, EB, je dis que comme le rectangle FAG est au rectangle FCG, ainsi le rectangle BAE est au rectangle BCE~». Effectuons la traduction en termes aguésiens. Les points $K,N,V,O$ forment un quadrilatère dont ils sont les \textit{bornes}, et les deux diagonales $KN,VO$ sont des \textit{bornales couplées,} de même que $VN$ et $KO$. Prenant une droite quelconque $\Delta$, celle-ci coupe les bornales en des points accouplés $A,C$ et $E,B$, ainsi que la conique contenant $K,N,V,O$ en $F,G$. On retrouve la combinatoire choisie par Beaugrand un peu plus haut dans  sa lettre, et il affirme qu'alors
\[
\frac{AF.AG}{CF.CG}=\frac{AB.AE}{CB.CE},
\]
ou, pour reprendre les termes de Desargues, que  les rectangles de couples de brins sont à leur relatif comme leurs gémeaux sont entre eux.

Nous n'entrerons pas ici dans le détail de la preuve donnée par Beaugrand, qui s'appuie, écrit-il, sur la proposition 17 du livre III des \textit{Coniques} d'Apollonius (\textit{voir} \cite{apollonius-rashed}, p. 312). Cette proposition traite de proportions faisant intervenir des \textit{tangentes} à la conique, et le moins que l'on puisse dire est que l'application de celle-ci par Beaugrand est quelque peu expéditive. 

Mais l'important pour nous est plus loin, aux lignes 10 et suivantes de la page 6 de l'original. Il y écrit qu'il peut de même démontrer que
\[
\frac{FA.FC}{GA.GC}=\frac{FB.FE}{GB.GE}.
\]
Il poursuit ainsi~:~«~L.S.D.~a aussi passé sous silence cette Analogie qui est plus considerable que les precedentes (et) il aura eu sans doute quelques \textit{evenemes de contrarietez,} en recherchant sa demonstration, bien qu'elle ne soit pas difficile~». Beaugrand n'a pas vu que cette identité n'est qu'une conséquence de l'énoncé de Desargues selon lequel les points $(AC,BE,CG)$ sont en involution, mais il nous révèle surtout que ce qui intéresse Beaugrand, ce qui pour lui \textit{fait géométrie,} ce sont les analogies, les égalités de rapports. Cette position épistémologique est incompatible avec celle de Desargues, qui s'intéresse lui avant tout aux \textit{dispositions de points.} Pour Desargues, l'important est la notion d'involution, disposition particulière de points sur une droite qui a son intérêt pour plusieurs raisons~:~
\begin{itemize}
\item elle s'exprime par des égalités de rapports, permettant l'utilisation de l'attirail euclidien;
\item elle intervient de manière naturelle et universelle dans la théorie des coniques;
\item elle est invariante par perspective;
\item elle est «~agrégative~» au sens défini plus haut.
\end{itemize}
En \emph{nommant} une configuration de points, comme celle d'involution, qui, par ses propriétés «~plastiques~», se rapproche de ce que deviendra l'involution considérée comme \emph{relation,}  Desargues opère une rupture que Beaugrand n'a pas vu ou peut-être pas pu accepter. Il écrit par exemple, p. 4 de sa lettre que «~c'est en ces trois Analogies que consiste ce que le S.D. nomme hors de propos, \textit{involution de six poincts,}~» ou encore, un peu plus bas~:~«~il y a une infinité de differentes analogies qui se peuvent former de quatre rectangles ou d'autres quantitez, si on vouloit introduire des termes pour les exprimer chacune en particulier il n'y en auroit pas assez en toutes les langues du monde pour en fournir autant qu'il seroit necessaire~». 

\begin{center}
{\Large\bf Conclusion}
\end{center}

\vs

La lecture détaillée des dix premières pages du \textit{Brouillon-Project} permet de comprendre en quoi l'approche de Desargues de la géométrie est novatrice. En ne se contentant pas seulement de démontrer des égalités de rapports mais en introduisant une notion nouvelle concernant la disposition de points sur une droite, notion qui enrichit et précise celle de division harmonique, il ouvre des possibilités nouvelles qui lui permettront, dans la seconde moitié du texte, de développer de manière unifiée la théorie des coniques. Si la riche terminologie que Desargues introduit n'est pas nécessaire et peut à juste titre être critiquée, la notion d'involution est cependant d'une radicale nouveauté, de par son agrégativité, que Beaugrand n'a pas comprise, mais aussi de par son invariance par projection centrale. Son long développement sur les n{\oe}uds moyens peut sembler indigeste mais la clarification qu'il apporte à la notion de division harmonique, et donc aux énoncés d'Apollonius, le rend lui aussi crucial.

De fait, Desargues va utiliser ces propriétés dans son extension et son unification de la théorie apollonienne des diamètres et des ordonnées. Cela est particulièrement saillant entre le milieu de la page 20 et celui de la page 22 du {\it Brouillon--project,} où il développe sa théorie des traversales\footnote{Que l'on peut apparenter aux polaires.}. Il y démontre en particulier que si trois points sont alignés, alors leurs trois traversales sont concourantes, ainsi que le résultat dual, ce qui s'exprime en termes modernes en disant que la polarité entre le plan projectif et son dual est une homographie, énoncé aujourd'hui à la base de l'étude des coniques. 

Il utilise {\it pleinement} l'agrégativité de l'involution pour arriver à ces propositions\footnote{{\it Voir} la fin de la page 20 et le début de la page 21 du \textit{Brouillon,} par exemple.}, autrement dit la phrase de synthèse de la page 10, que nous citons à nouveau~:~«~Quand en un tronc droict, trois couples de n{\oe}uds extremes DF, CG, BH, sont en involution entre eux, que deux autres couples de n{\oe}uds, moyens, unis, doubles ou simples, PQ,XY, sont une involution de quatre points avec chacune des deux quelconques couples CG,\& BH, de ces trois couples de n{\oe}uds extremes. Ces deux mesmes n{\oe}uds moyens PQ, XY, sont encore une involution aussi de quatre points avec la troisiéme de ces couples de n{\oe}uds extremes DF.~»

Mais il va encore plus loin en écrivant, aux lignes 22 et suivantes de la page 22 du {\it Brouillon--Project~:~}«~D'où suit aussi qu'autant de couples de droictes qui sont ordonnées à un des poincts du bord de la coupe de rouleau, \& qui passent aux deux poincts du mesme bord qu'y donne une quelconque droicte d'une quelconque ordonnance, donnent en la traversale de cette ordonnance {\it autant de couples de n{\oe}uds d'une involution}.~» On peut être tenté de voir ici en germe la notion de \textit{correspondance} ou de transformation, mais cette notion reste évanescente, elle se dérobe aux yeux de Desargues qui ne peut donc totalement clarifier ses conceptions de l'involution. Cette impression que Desargues s'approche du concept de transformation, sans toutefois arriver à le saisir, est renforcée par son insistance parfois lourde sur la \textit{combinatoire} de l'involution, combinatoire qui contient en germe la notion de transformation involutive (au travers des couples \textit{relatifs}) et qui ne se réduit pas, chez Desargues, à un simple outil technique de génération de nouvelles analogies. Nous montrerons dans un article à venir ({\it voir} \cite{anglade-Briend-2}) comment cette emploi de la combinatoire structure la pratique mathématique de Desargues en analysant son usage du théorème de Ménélaüs. 

Ce qui en tout cas ne fait aucun doute, c'est que les conceptions unificatrices concernant les coniques, plus particulièrement son développement d'idées géométriques qu'il faut bien qualifier de \textit{projectives}, ont profondément influencé Blaise Pascal pour son \textit{Essay sur les coniques,} et que ce dernier a sans doute été l'élève direct du lyonnais.

L'on peut donc en conclure que la position épistémologique de Desargues est en rupture avec celle de son temps concernant la géométrie~:~d'une conception centrée sur les identités de rapports, on passe à une conception centrée sur des configurations qui possèdent des propriétés intéressantes d'invariance par transformation, comme une projection centrale par exemple. Le côté novateur de cette approche, son aspect inachevé, la lourdeur du style de Desargues et le déclin de l'intérêt pour la théorie des coniques peuvent expliquer en partie pourquoi l'{\oe}uvre de Desargues a connu une si longue éclipse après son surgissement sur la scène mathématique.

\vspace*{1cm}

\noindent{\large \textbf{Appendice~:~les résultats euclidiens utilisés par Desargues pour son étude de l'involution}}

\vs
Dans les démonstrations de ses résultats fondamentaux sur l'involution, à savoir l'équivalence entre les notions d'arbre et d'involution ({\it voir.} le paragraphe 2.3) et le théorème de la Ramée énonçant l'invariance par projection centrale, Desargues utilise un langage et des méthodes strictement euclidiens. Il va employer sans jamais les nommer, car ils font partie du bagage commun à tout lecteur potentiel du {\it Brouillon,} quelques résultats de bases sur les proportions que nous rappelons ici pour la commodité du lecteur.

En revanche, à la fin de la page 2, il donne une liste explicite de propositions tirées des {\it Éléments} d'Euclide, qui curieusement, ne seront nullement utilisées dans les dix premières pages. Donnons-en la citation complète~:«~Proposition comprenant les 5 \& 6 du second des Elemens d'Euclides. Proposition comprenant les 9 \& 10 du 2 des Elemens d'Euclides. Proposition comprenant les 35 \& 36 du 3 des Elemens d'Euclides~». Desargues précise dans ses notes correctives qu'il s'agit de ces propositions et de leur {\it converse,} c'est-à-dire leur réciproque. Il va en fait les rappeler à la page 26 du {\it Brouillon} pour en faire usage dans son étude des propriétés des asymptotes de l'hyperbole. Nous les rappellerons cependant, afin de donner un exposé complet. Il cite enfin\footnote{p. 2, à partir de la dernière ligne.} une proposition, «~énoncée autrement en Ptolomée\footnote{p. 3, l. 4}~» qui n'est autre que le théorème aujourd'hui connu sous le nom de {\it théorème de Ménélaüs} et dont il donne une preuve un peu plus loin. Nous laisserons de côté cette dernière proposition, qui ne trouve son usage que lorsque Desargues s'attaque à la géométrie plane proprement dite, lors de la démonstration du théorème de la Ramée énonçant l'invariance de la configuration d'involution par projection centrale.

\vs
Il est difficile de savoir à quelle source Desargues se réfère concernant les résultats euclidiens, et nous nous sommes basés sur l'édition française des {\it Éléments} parue en 1621, due à Didier (ou Denis) Henrion ({\it voir} \cite{henrion}). Comme à son habitude, il évite la terminologie latine ou d'origine grecque et utilise à la place un vocabulaire français. On trouve souvent à la fin d'un argument concernant des proportions une conclusion comme «~\& alternement, composant, divisant, \& le reste~». Cela signifie qu'en utilisant les propositions classiques d'Euclide concernant les proportions, souvent cités sous leur nom latin de {\it invertendo, componendo} ou {\it dividendo,} on peut déduire un certain nombre de résultats qu'il pourra utiliser par la suite.

\vs
\noindent{\it Alternando.~} Étant donnée la proportion $a/b=c/d$, alors $a/c=b/d$.

\noindent{\it Componendo.~} Étant donnée la proportion $a/b=c/d$, alors $(a+b)/b=(c+d)/d$.

\noindent{\it Dividendo.~} Étant donnée la proportion $a/b=c/d$, alors $(a-b)/b=(c-d)/d$.

\noindent{\it Invertendo.~} Étant donnée la proportion $a/b=c/d$, alors $b/a=d/c$.

\noindent{\it Convertendo.~} Étant donnée la proportion $a/b=c/d$, alors $a/(a-b)=c/(c-d)$.

\vs
\noindent Terminons cette courte liste par un résultat que Desargues utilise souvent, qui est la proposition 12 du livre 5 des {\it Éléments,} que l'on peut énoncer ainsi~:~si
\[
\frac{a}{b}=\frac{c}{d}
\]
alors
\[
\frac{a}{b}=\frac{a+c}{b+d},
\]
ainsi que la proportion qui s'en déduit par soustraction, et alternement, composant, divisant, \& le reste.

\bibliographystyle{plain}

\newpage

\noindent \textbf{Remerciements~:}~les auteurs tiennent à remercier l'équipe de la licence \textit{sciences et humanités} de l'université d'Aix-Marseille, sans qui ce travail n'aurait jamais vu le jour ; nous remercions en particulier Sara Ploquin-Donzenac pour son aide précieuse et constante. Nous tenons également à remercier  Valérie Debuiche et Sylvie Pic pour les fructueuses discussions que nous avons au sujet du \textit{Brouillon Project}. Enfin, nous remercions plus particulièrement Philippe Abgrall pour ses nombreuses suggestions et sa relecture attentive du présent texte.

\vs
\noindent Marie Anglade\\ 
Université d'Aix-Marseille\\
CEPERC UMR CNRS 7304\\
3, place victor Hugo\\
13331 Marseille cede 3\\
France. \\
\verb+marie.anglade@univ-amu.fr+

\vs
\noindent Jean-Yves Briend\\
Université d'Aix-Marseille\\
I2M UMR CNRS 7373\\
CMI\\
39, rue Joliot-Curie\\
13453 Marseille cedex 13\\
France.\\
\verb+jean-yves.briend@univ-amu.fr+

\end{document}